# A Finite Volume Method for Elastic Waves in Heterogeneous, Anisotropic and Fractured Media

Ingrid Kristine Jacobsen[1], Inga Berre[1], Jan Martin Nordbotten[1,2], Ivar Stefansson[1]

[1] Center for Modeling of Coupled Subsurface, Dynamics, Department of Mathematics, University of Bergen, Allégaten 41, 5020 Bergen, Norway.

[2] NORCE Norwegian Research Centre AS

## Abstract

Numerical modeling of elastic wave propagation in the subsurface requires applicability to heterogeneous, anisotropic and discontinuous media, as well as support of free surface boundary conditions. Here we study the cell-centered finite volume method Multi-Point Stress Approximation with weak symmetry (MPSA-W) for solving the elastic wave equation. Finite volume methods are geometrically flexible, locally conserving and they are suitable for handling material discontinuities and anisotropies. For discretization in time we have utilized the Newmark method, thereby developing an MPSA-Newmark discretization for the elastic wave equation. An important aspect of this work is the integration of absorbing boundary conditions into the MPSA-Newmark method to limit possible boundary reflections.

We verify the MPSA-Newmark discretization numerically for model problems. Convergence analysis of MPSA-Newmark is performed using a known solution in a medium with homogeneous Dirichlet boundary conditions. The analysis demonstrates the expected convergence rates of second order for primary variables (displacements) and between first and second order for secondary variables (tractions). Further verification is conducted through convergence analysis with the inclusion of absorbing boundary conditions. The stability of the scheme is shown through numerical energy decay analyses for waves travelling with various incidence angles onto the absorbing boundaries. Lastly, we present simulation examples of wave propagation in fractured, heterogeneous and transversely isotropic media to demonstrate the versatility of the MPSA-Newmark discretization.

Data availability: The source code for all simulations is open-source and it is available at Zenodo (Jacobsen et al., 2025).

## 1 Introduction

Modeling of elastic waves in the subsurface is important in understanding natural and induced earthquakes. Induced seismicity, e.g. related to geothermal energy extraction, $CO_2$-storage and wastewater disposal, takes place in locations characterized by complex material compositions. Specifically, heterogeneous and anisotropic rock formations are highly relevant in all the abovementioned fields of application. Solving elastic wave propagation within such complex rock formations means that the numerical method must handle unstructured and general grids and material anisotropy and heterogeneity. Moreover, simulations involving the Earth's surface require employment of the free surface boundary condition.

Several numerical methods have been developed and utilized to solve the elastic wave propagation problem. Commonly studied methods within the research literature include those within the families of spectral element (SEM) and discontinuous Galerkin (dG) methods. Other methods such as finite difference (FD), pseudospectral (PSM), finite elements (FEM) and finite volume (FV) methods have also been studied. All these methods have their own advantages and drawbacks. For instance, finite difference and spectral element methods are very efficient for elastic wave propagation problems, but they perform best in media without strong structural complexities (e.g. fractures) and material heterogeneities (Igel, 2016). On the other hand, finite volume and discontinuous Galerkin methods are naturally applicable to discontinuous and heterogeneous media (Cardiff & Demirdžić, 2021; Keilegavlen & Nordbotten, 2017; Virieux et al., 2011), and they are flexible in terms of arbitrary cell shapes and unstructured grids (Cardiff & Demirdžić, 2021). For an overview of numerical methods for the elastic wave equation, we refer to Igel (2016) and Seriani and Oliveira (2020).

Finite volume and discontinuous Galerkin methods are closely related and share many advantageous properties. As already mentioned, they support a wide range of grids, including unstructured ones, which allow for representation of complex simulation domains. Such grid flexibility is crucial for accurately representing fractures and heterogeneous media where the material discontinuities do not coincide with orthogonal gridlines, as well as other geometrical complexities such as surface topography or boreholes (Igel, 2016). In addition to this, finite volume and discontinuous Galerkin methods can naturally enforce the free surface boundary condition, which is of high relevance in simulations including the Earth's surface. Although finite volume methods offer some advantageous properties, they are not among the most widely used for solving elastic wave propagation problems. This is due to finite volume methods often being restricted to lower order, while discontinuous Galerkin methods are easily extended to high order accuracy (Igel, 2016). However, the higher order accuracy comes at the expense of additional computational cost, which may not be appropriate when geometric features reduce the regularity of the solution.

As finite volume methods are widely used within the realm of flow in porous media, developing and utilizing finite volume methods for mechanical deformation makes it possible to have a unified solution approach for poroelastic solid-fluid interaction problems (Nordbotten, 2014). Common solution procedures for such coupled problems often include solving each of the subproblems with separate discretization schemes. See e.g. Kim et al. (2011) for an example of a finite element method for the mechanics and a finite volume method for the flow. Several works have focused on utilizing a finite volume method for both flow and mechanical deformation, providing a unified solution approach and consistent data structure for both subproblems (Nordbotten, 2016; Novikov et al., 2022; Terekhov, 2020; Terekhov & Vassilevski, 2022). All these works consider cell-centered finite volume methods for the coupled flow and static mechanical deformation problem. A natural extension is thus to develop finite volume methods for dynamic mechanical deformation/elastic wave propagation which can later, in a fully coupled manner, be used to solve poroelastic solid-fluid interaction problems with the inclusion of seismic effects.

There are several finite volume methods for mechanical deformation, and as we will highlight only a few, we refer to Cardiff and Demirdžić (2021) for an extensive review. The first application of a cell-centered finite volume method to solid mechanics was a method for homogeneous and isotropic linear elasticity by Demirdzic et al. (1988). Early developments of finite volume methods for elastic wave propagation include that of Dormy and Tarantola (1995) which use the velocity-displacement formulation of the elastic wave equation. Tadi (2004) proposed a second order accurate method for the second order displacement formulation with a focus on natural application of traction boundary conditions. Higher order finite volume methods for unstructured triangular and tetrahedral meshes include the arbitrary high order method by Dumbser et al. (2007). A more recent development is that of Zhang et al. (2016), which combines techniques of high order and spectral finite volume methods for the 2D problem. Zhang et al. (2017) presented an extension of this methodology to three spatial dimensions. Lemoine et al. (2013) used a high-resolution finite volume method for wave propagation in orthotropic poroelastic media on cartesian grids, which was later expanded to nonrectilinear mapped grids by Lemoine and Ou (2014). Another finite volume method used for static mechanical deformation is the work of Tuković et al. (2013), which presents a detailed description of the finite volume discretization of multiple deformable bodies with different solid material parameters. Cell-centered finite volume methods for anisotropic and heterogeneous media include the multi-point stress approximation with weak symmetry (MPSA-W) by Keilegavlen and Nordbotten (2017) and the method of Terekhov and Tchelepi (2020).

The finite volume MPSA-W method, which is up to second order accurate, applies to a broad range of grids, can naturally enforce the free surface boundary condition and supports heterogeneous and anisotropic media. In addition to this, MPSA-W is the vector problem analogue of the finite volume Multi-Point Flux Approximation (MPFA) (Nordbotten & Keilegavlen, 2021; Aavatsmark, 2002), a spatial discretization that is shown to be a successful discretization method for solving problems within flow in porous media (Berre et al., 2021). MPSA-W is built on the same framework as MPFA, meaning the two methods hold the same advantages of computational efficiency due to minimal number of degrees of freedom and being able to solve problems in discontinuous, heterogeneous and anisotropic media. Thus, the MPFA and the MPSA-W are examples of methods that can be combined for a seamlessly coupled solution procedure for poroelastic fluid-solid interaction problems. The MPFA-MPSA pair is available in simulation toolboxes such as the Porous Microstructure Analysis (PuMA) (Ferguson et al., 2018), the MATLAB Reservoir Simulation Toolbox (MRST) (Lie, 2019), as well as in PorePy (Keilegavlen et al., 2021) which is used herein.

Motivated by the properties of the MPSA-W method for (quasi)static elastic and poroelastic deformation problems in heterogeneous and anisotropic media, we present here for the first time a study of the MPSA-W method for solving the elastic wave propagation problem. The spatial discretization is applied in combination with a Newmark time integration scheme (Newmark, 1959), which is a second order accurate implicit time integration scheme widely used in computational solid mechanics. By combining the spatial discretization, MPSA-W, and the temporal discretization, Newmark, we obtain the MPSA-Newmark method. The MPSA-Newmark method is up to second order accurate in space and time, meaning it has balanced space/time accuracy for solving the elastic wave equation.

Solving the elastic wave equation poses the question of how to deal with possible wave reflections on the domain boundary. We have addressed this by employing the first order absorbing boundary conditions (ABCs) first presented by Clayton and Engquist (1977), which in practice are time-dependent Robin-type boundary conditions. Therefore, this paper also provides a presentation of how time-dependent Robin boundary conditions are discretized in MPSA. As the structure of MPSA is like that of MPFA, the treatment of the boundary conditions is analogous in both discretization schemes. This analogy extends the applicability of the presented methodology for discretizing time-dependent Robin conditions to other multi-point finite volume schemes.

The structure of the paper is as follows. In Section 2, we present the mathematical model and the discretization. In Section 3 we present verification of the code through convergence analyses for two different model problems: one regarding homogeneous Dirichlet boundary conditions, and one regarding the inclusion of an absorbing boundary. Section 3 also presents a numerical analysis of the scheme stability through energy decay investigations in an isotropic media with all absorbing boundaries. In Section 4 we show simulation examples of wave propagation in anisotropic, heterogeneous and fractured media. Concluding remarks are given in Section 5.

## 2 Methodology

This section covers the mathematical model and the discretization of the problem. We refer to Table 2.1 for an overview of the symbols used in the article.

*Table 2.1: Nomenclature.*

| Symbol | Physical parameter | Unit |
|---|---|---|
| $n$ | Unit normal vector | – |
| $v$ | Symmetry axis vector for transversely isotropic media | |
| $\mathcal{C}$ | Stiffness tensor | kg $\mathrm{m}^{-1}\mathrm{s}^{-2}$ |
| $\mathcal{R}$ | Robin weight | – |
| $\mathcal{D}, \boldsymbol{\mathcal{D}}$ | Robin weight for absorbing boundary conditions, discrete Robin weight for absorbing boundary conditions | kg $\mathrm{m}^{-2}\mathrm{s}^{-1}$ |
| $\mathcal{D}_{\Delta t}, \boldsymbol{\mathcal{D}}_{\Delta t}$ | $\mathcal{D}$ scaled with inverse of the time-step, discrete $\boldsymbol{\mathcal{D}}$ scaled with inverse of the time-step | kg $\mathrm{m}^{-2}\mathrm{s}^{-2}$ |
| $q, \boldsymbol{q}$ | Body force, discrete body force | kg $\mathrm{m}^{-2}\mathrm{s}^{-2}$ |
| $T, \boldsymbol{T}$ | Traction, discrete traction | kg $\mathrm{m}^{-1}\mathrm{s}^{-2}$ |
| $u, \boldsymbol{u}$ | Displacement, discrete displacement | m |
| $\dot{u}, \dot{\boldsymbol{u}}$ | Velocity, discrete velocity | m $\mathrm{s}^{-1}$ |
| $\ddot{u}, \ddot{\boldsymbol{u}}$ | Acceleration, discrete acceleration | m $\mathrm{s}^{-2}$ |
| $\sigma$ | Stress tensor | kg $\mathrm{m}^{-1}\mathrm{s}^{-2}$ |
| $\epsilon$ | Symmetric gradient of the displacement | m $\mathrm{m}^{-1}$ |
| $\tau$ | Internal traction weights | kg $\mathrm{m}^{-2}\mathrm{s}^{-2}$ |
| $w$ | Boundary traction weights | kg $\mathrm{m}^{-2}\mathrm{s}^{-2}$ |
| $F$ | Right-hand side for the boundary conditions (with subscript to denote boundary condition type) | m or kg $\mathrm{m}^{-1}\mathrm{s}^{-2}$ |
| $\rho$ | Density | kg $\mathrm{m}^{-3}$ |
| $\lambda$ | First Lamé parameter | kg $\mathrm{m}^{-1}\mathrm{s}^{-2}$ |
| $\lambda_{\|\|}$ | Transverse compressive stress parameter | kg $\mathrm{m}^{-1}\mathrm{s}^{-2}$ |
| $\lambda_{\perp}$ | Perpendicular compressive stress parameter | kg $\mathrm{m}^{-1}$ $\mathrm{s}^{-2}$ |
| $\mu$ | Second Lamé parameter/Shear modulus | kg $\mathrm{m}^{-1}\mathrm{s}^{-2}$ |
| $\mu_{\|\|}$ | Transverse shear parameter | kg $\mathrm{m}^{-1}\mathrm{s}^{-2}$ |
| $\mu_{\perp}$ | Transverse-to-perpendicular shear parameter | kg $\mathrm{m}^{-1}$ $\mathrm{s}^{-2}$ |
| $c$ | Primary wave speed | m $\mathrm{s}^{-1}$ |
| $\theta$ | Wave rotation angle | rad |
| $\mathcal{E}_u$ | Relative discrete $L^2$-error of displacement | – |
| $\mathcal{E}_T$ | Relative discrete $L^2$-error of traction | – |
| $E$ | Energy | kg m$^2$s$^{-2}$ |
| $E_0$ | Initial energy | kg m$^2$s$^{-2}$ |
| $\Omega$ | Simulation domain | m$^2$ or m$^3$ |
| $\omega$ | Anisotropic region in $\Omega$ | m$^2$ or m$^3$ |
| $L$ | Width of domain A in convergence analysis | m |
| $r_h$ | Heterogeneity coefficient | – |
| $r_a$ | Anisotropy coefficient | – |
| $t$ | Time | s |
| $\beta$ | Time-discretization parameter | – |
| $\gamma$ | Time-discretization parameter | – |
| $\Delta t$ | Time-step size | s |
| $n$ (as superscript) | Time-step number | – |
| $\Delta x$ | Grid size | m |
| $N_x$ | Number of cells | – |
| $N_t$ | Number of time-steps | – |
| $\nabla$ | Nabla-operator | m$^{-1}$ |

## 2.1 Governing Equations

### 2.1.1 The Elastic Wave Equation

We consider propagation of elastic waves in a three-dimensional domain denoted $\Omega$, the reduction to two spatial dimensions being straightforward. As our interest is in elastic waves, we will only consider infinitesimal deformations in an elastic material. The equation for conservation of momentum then takes the form:

$$\rho \ddot{u} = \nabla \cdot \sigma + q \qquad \text{in } \Omega. \tag{2.1}$$

Here $\rho$ is the rock density, $u = [u_x, u_y, u_z]^{\mathrm{T}}$ is the displacement, $\sigma$ is the stress tensor, $q$ is an external body force and, using the dot notation for derivatives in time, $\ddot{u}$ is the acceleration. The constitutive stress-strain relationship is given for small deformations by Hooke's law as:

$$\sigma = C : \epsilon(u), \tag{2.2}$$

where $C$ is the fourth order stiffness tensor and $\epsilon(u)$ is the symmetric gradient of $u$:

$$\epsilon(u) = \frac{1}{2}(\nabla u + (\nabla u)^T). \tag{2.3}$$

For future reference, we present two examples of the stiffness tensor $C$. First, we present $C$ for an isotropic medium:

$$C_{ijkl} = \lambda \delta_{ij}\delta_{kl} + \mu(\delta_{ik}\delta_{jl} + \delta_{il}\delta_{jk}), \tag{2.4}$$

where $\lambda$ is the first Lamé parameter, $\mu$ is the second Lamé parameter, or the shear modulus, and $\delta_{ij}$ is the Kronecker delta, taking the value 1 if $i = j$, and 0 otherwise. Secondly, we will present $C$ for a transversely isotropic medium (Payton, 2012). Transversely isotropic media are represented by five independent material parameters, as opposed to only two as is the case for isotropic media. The stiffness tensor for transversely isotropic media, where the axis of symmetry aligns with the unit vector $v = (v_1, v_2, v_3)^T$, is as follows:

$$
\begin{aligned}
\mathcal{C}_{ijkl} = (\lambda + \lambda_{||})\delta_{ij}\delta_{kl} + \mu_{||}(\delta_{ik}\delta_{jl} + \delta_{il}\delta_{jk}) - \lambda_{||}(\delta_{ij}v_k v_l + \delta_{kl}v_i v_j) \\
+ (\mu_{\perp} - \mu_{||})(\delta_{ik}v_j v_l + \delta_{jk}v_i v_l + \delta_{il}v_j v_k + \delta_{jl}v_i v_k) \\
+ (\lambda_{||} + \lambda_{\perp} + 2\mu_{||} - 2\mu_{\perp})v_i v_j v_k v_l, \quad \text{for} \quad i,j,k,l = 1,2,3.
\end{aligned} \tag{2.5}
$$

Here $\lambda$ is the first Lamé parameter, $\lambda_{||}$ and $\lambda_{\perp}$ are the transverse and perpendicular compressive stress parameters, respectively, and $\mu_{||}$ and $\mu_{\perp}$ are the transverse and transverse-to-perpendicular shear parameters, respectively. Notice that by setting $\lambda_{||} = \lambda_{\perp} = 0$ and $\mu_{||} = \mu_{\perp}$, the dependence on $v$ is eliminated from Equation (2.5). In that case, Equation (2.5) reduces to Equation (2.4), that is, the stiffness tensor for an isotropic solid.

To close the mathematical model, we have initial and boundary conditions. The boundary conditions are detailed in dedicated sections below, while the initial conditions are included in Section 2.2.3.

### 2.1.2 Boundary conditions

We denote the boundary of the domain as $\Gamma$, which we subdivide into four non-overlapping parts, $\Gamma = \Gamma_N \cup \Gamma_D \cup \Gamma_R \cup \Gamma_A$. These are identified as the Neumann, Dirichlet, Robin and Absorbing boundaries, respectively. We summarize these below.

Neumann boundary conditions are specified in terms of prescribed tractions, $F_N$:

$$
\sigma \cdot n = F_N \qquad \text{on } \Gamma_N, \tag{2.6a}
$$

where $n$ is the outward pointing normal vector relative to the boundary. Dirichlet boundary conditions are specified in terms of prescribed displacements, $F_D$:

$$
u = F_D \qquad \text{on } \Gamma_D. \tag{2.6b}
$$

Robin conditions are specified as a weighted sum of traction and displacement at the boundary, with prescribed data $F_R$ having units of traction:

$$
\sigma \cdot n + \mathcal{R}u = F_R \qquad \text{on } \Gamma_R. \tag{2.6c}
$$

The coefficient $\mathcal{R}$ is a 2-tensor weight for the displacement component of the Robin conditions.

It is often desirable to eliminate wave reflections at the boundaries when solving the elastic wave equation. A natural choice is to employ boundary conditions that effectively allow outward-going waves to pass through without reflection, often referred to as absorbing boundary conditions. We consider here one of the simplest forms of absorbing boundary conditions, namely the first order ones by Clayton and Engquist (1977) which absorb waves at normal incidence exactly (Higdon, 1991). The first order absorbing boundaries cannot perfectly absorb waves incident at non-orthogonal angles. However, they do reduce the amplitude of reflected waves from such angles, unlike Dirichlet and Neumann conditions, which fully reflect the waves. The absorbing conditions by Clayton and Engquist were originally expressed with displacements and displacement gradients, but were later expressed in terms of displacements and tractions by Tsogka (1999). We utilize the latter formulation

$$
\sigma \cdot n + \mathcal{D}\dot{u} = F_A \qquad \text{on } \Gamma_A, \tag{2.7}
$$

where $\mathcal{D} = \left(\sqrt{\rho(\lambda + 2\mu)}\, nn^T + \sqrt{\rho\mu}(I - nn^T)\right)$. The time-dependent component of the boundary condition is the time derivative of the boundary displacement, namely the velocity, $\dot{u}$. For the sake of generality, we have included a non-zero forcing term $F_A$ with the understanding that the absorbing boundary condition is recovered by setting $F_A = 0$.

When referring to the boundary data in general, we will sometimes simply refer to $F$, understanding that the subscript should be clear from context.

## 2.2 Discretization

In this section we will summarize the main structure of the space-time discretization, emphasizing the novel aspects of this work.

### 2.2.1 Multi-Point Stress Approximations Including Boundary Conditions

To simplify the exposition, we restrict the presentation to polyhedral domains $\Omega$. We then consider a non-overlapping set $\mathcal{T}$ of $k$ cells $K_i \in \mathcal{T}$ for $i = 1, 2, \dots, k$, such that $\bigcup_{K_i \in \mathcal{T}} \overline{K_i} = \Omega$. A face on the boundary of a cell is denoted by $f$ and is collected in the set of all faces, $\mathcal{F}$. The set of faces for a particular cell $K$ is denoted $\mathcal{F}_K$. In a similar manner, neighboring cells of a face $f$ are denoted $\mathcal{T}_f$. Figure 2.1 provides an illustration of an internal face $f_1$ (red) and the corresponding $\mathcal{T}_{f_1}$ (marked in pink). Cell $K$ has volume $m_K$, and face $f$ has measure $m_f$. Boundary related grid quantities follow the same notation as mentioned above, but with a tilde. For instance, $\tilde{f}$ represents a boundary face, and the set of all boundary faces is $\tilde{\mathcal{F}} \subset \mathcal{F}$. The set of boundary faces sharing a vertex with a face $f$ is denoted $\tilde{\mathcal{F}}_f$. We refer to Figure 2.1 for an illustration of the face $f_2$ (blue) which touches the boundary and the corresponding $\tilde{\mathcal{F}}_{f_2}$ (marked in black). Finally, sets of Dirichlet, Neumann, Robin or Absorbing boundary faces will have subscripts $D$, $N$, $R$ or $A$, analogously to their continuous counterparts. Figure 2.1 illustrates a two-dimensional unstructured grid with the grid-related quantities we will use to present MPSA and boundary conditions, but the same relations hold also in three spatial dimensions.

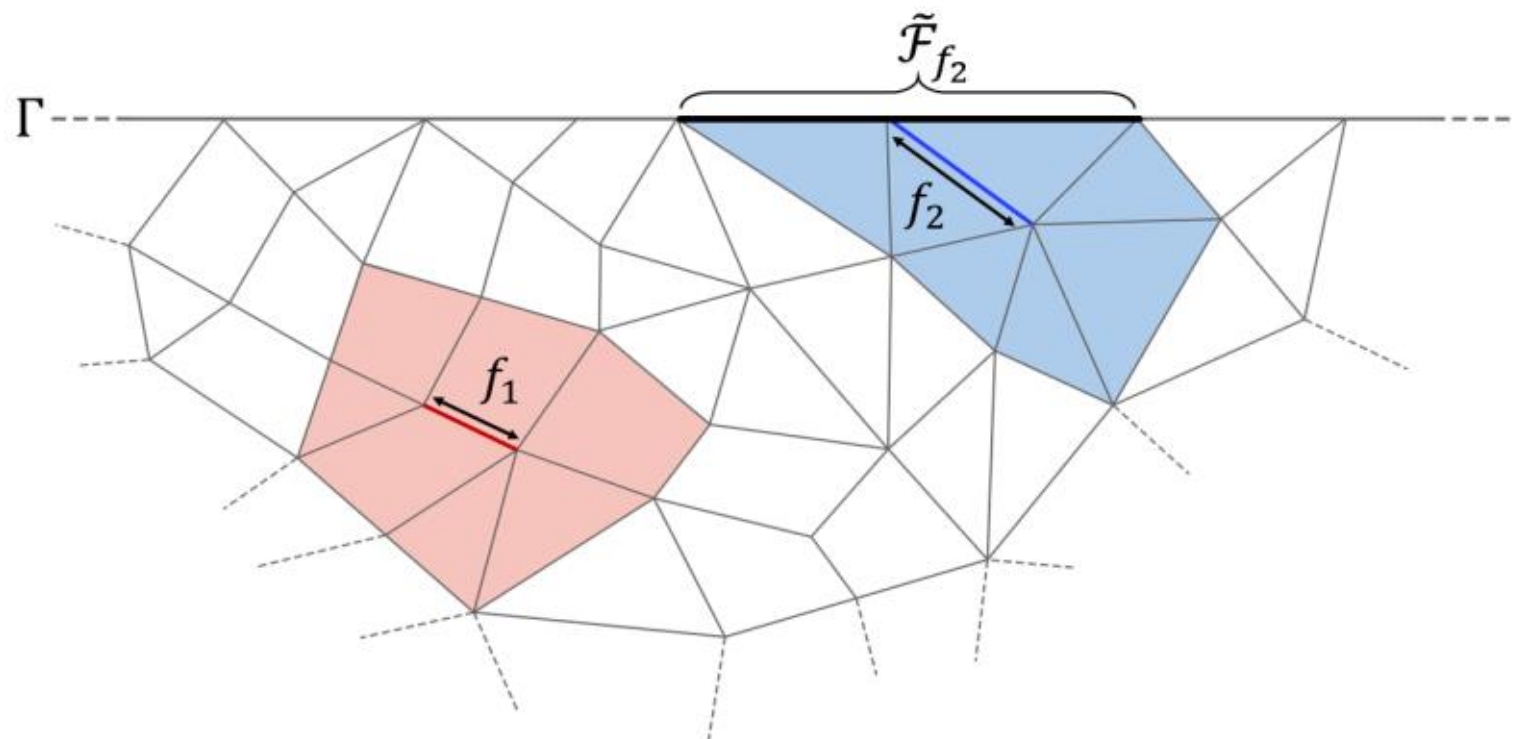

*Figure 2.1: Illustration of grid-related quantities used to express the Multi-Point Stress Approximation and the discrete boundary conditions. The pink shaded area is the stencil for approximating the flux over the internal face $f_1$ (marked in red). The blue shaded area and the boundary faces $\tilde{\mathcal{F}}_{f_2}$ (marked in black) represent the stencil for approximating the flux over the face $f_2$ (marked in blue) that has an edge touching the boundary $\Gamma$.*

Integrating Equation (2.1) over a cell $K \in \mathcal{T}$ and applying the divergence theorem gives

$$\int_K (\rho \ddot{u} - q)\, dV = \int_{\partial K} \sigma \cdot n_{\partial K}\, dA, \tag{2.8}$$

where $n_{\partial K}$ is the outward pointing normal vector of the boundary $\partial K$ of cell $K$.

Dividing the surface integral into two sums of face tractions leads to the following general expression which holds for any finite volume method:

$$\int_{\partial K} \sigma \cdot n_{\partial K}\, dA = \sum_{f \in \mathcal{F}_K \setminus \tilde{\mathcal{F}}_{K,N}} \int_f \sigma \cdot n_{K,f}\, dA + \sum_{\tilde{f} \in \tilde{\mathcal{F}}_{K,N}} \int_{\tilde{f}} \sigma \cdot n_{K,\tilde{f}}\, dA. \tag{2.9}$$

Here we have separated the sum into faces that lie on a Neumann boundary (where the traction is known directly from the boundary conditions) and remaining faces. We now define cell-averaged displacements and face-averaged tractions as:

$$\boldsymbol{u}_K = \frac{1}{m_K} \int_K u\, dV \qquad \text{and} \qquad \boldsymbol{T}_{K,f} = \frac{1}{m_f} \int_f \sigma \cdot n_{K,f}\, dA.$$

Here the bold font indicates that $\boldsymbol{u}$ and $\boldsymbol{T}$ are finite-dimensional vectors. Moreover, we note that for any two cells $K$ and $K'$ sharing a face $f$, the definition of the face traction only differs by the orientation of $n_{K,f}$. Therefore, due to the continuity of traction, it is clear that $\boldsymbol{T}_{K,f} = -\boldsymbol{T}_{K',f}$. Defining the cell-averaged source term analogously to the cell-averaged displacements, and considering that density is piece-wise constant, Equation (2.8) becomes the space-discrete (but still continuous in time) system of ordinary differential equations:

$$m_K(\rho \ddot{\boldsymbol{u}}_K - \boldsymbol{q}_K) = \sum_{f \in \mathcal{F}_K \setminus \tilde{\mathcal{F}}_{K,N}} m_f \boldsymbol{T}_{K,f} + \sum_{\tilde{f} \in \tilde{\mathcal{F}}_{K,N}} m_{\tilde{f}}\, \boldsymbol{F}_{\tilde{f}}, \tag{2.10}$$

where $\ddot{\boldsymbol{u}}_K$ and $\boldsymbol{q}_K$ are, respectively, the discrete acceleration and source for cell $K$, and $\boldsymbol{F}_{\tilde{f}}$ is the face-averaged Neumann boundary data for face $\tilde{f} \subset \Gamma_N$.

Equation (2.10) represents conservation of momentum in each finite volume $K$, and holds for all (spatial) finite volume approximations to the elastic wave equations. What distinguishes the various methods is how Hooke's law is approximated, in the sense of expressing $\boldsymbol{T}_{K,f}$ in terms of cell displacements $\boldsymbol{u}$. Most finite volume methods for elasticity, including the MPSA methods, are designed to be both linear as well as "local". The methods are local in the sense that, in the interior of the domain, the approximation of $\boldsymbol{T}_{K,f}$ only depends on the set $\mathcal{T}_f \subset \mathcal{T}$ of cells that share a vertex with the face $f$ (see the pink shaded region in Figure 2.1). For faces touching the boundary, like $f_2$ (the face colored in blue in Figure 2.1), the approximation of $\boldsymbol{T}_{K,f}$ includes the set $\tilde{\mathcal{F}}_f \subset \mathcal{F} \cap \Gamma$ of boundary faces sharing a vertex with $f$. The set $\tilde{\mathcal{F}}_f$ is empty for faces that do not share a node with the boundary. We limit the exposition to this class of linear and local methods, which implies that each traction can be written as:

$$\boldsymbol{T}_{K,f} = \sum_{L \in \mathcal{T}_f} \tau_{K,f,L} \boldsymbol{u}_L + \sum_{g \in \tilde{\mathcal{F}}_f} w_{K,f,g} \boldsymbol{F}_g. \tag{2.11}$$

The precise definition of the internal weights $\tau_{K,f,L}$ is not essential to the current presentation, and we refer to Keilegavlen and Nordbotten (2017) for details. We emphasize that the definition of the boundary weights $w_{K,f,g}$ depends on the type of boundary. Again, the previously mentioned reference is sufficient for defining the coefficients associated with the time-independent boundary conditions (i.e. Neumann, Dirichlet and Robin). However, that work did not consider absorbing boundary conditions, and we will detail the inclusion of such boundary conditions below.

The finite volume methods considered herein do not deal with the time-dependent aspect of boundary conditions, as they are exclusively defined as spatial discretizations. Therefore, the time derivative in Equation (2.7) must be discretized in time. Anticipating that we will use a second order accurate implicit scheme to discretize the acceleration term, we choose to consider a two-step implicit second order scheme for the boundary conditions, providing enhanced accuracy over a first-order method. This provides the following time-discrete representation of the absorbing boundary conditions:

$$\sigma^n \cdot n + \mathcal{D}_{\Delta t} u^n = \mathcal{D}_{\Delta t}\left(\frac{4}{3}u^{n-1} - \frac{1}{3}u^{n-2}\right) + F_A, \tag{2.12}$$

where $\mathcal{D}_{\Delta t} = \frac{3}{2\Delta t}\mathcal{D}$ and $\Delta t$ denotes the time-step size. The symbol $n$ is used for both the normal vector and the time-step number, where the distinction lies in whether $n$ is a superscript. When used in a superscript, $n$ denotes the time-step number. Quantities with superscript $(n-k)$ correspond to the quantity at $k$ time-steps back in time.

We recognize that the implicit problem to be solved, Equation (2.12), is on the form of a Robin boundary condition (see Equation (2.6c)) where the Robin weight $\mathcal{D}_{\Delta t}$ scales as the inverse of the time-step. As such, we note that the implementation of the absorbing boundary conditions is from a spatial perspective equivalent to implementing Robin boundary conditions with a non-zero data term obtained from the previous two time-steps. The fully discrete equation for the tractions on non-Neumann boundaries is thus:

$$\begin{aligned} \boldsymbol{T}_{K,f}^n = \sum_{L\in\mathcal{T}_f} \tau_{K,f,L}\boldsymbol{u}_L^n + \sum_{g\in\tilde{\mathcal{F}}_f\cap(\Gamma\setminus\Gamma_A)} w_{K,f,g}\boldsymbol{F}_g^n \\ + \sum_{g\in\tilde{\mathcal{F}}_f\cap\Gamma_A} w_{K,f,g}\left(\boldsymbol{\mathcal{D}}_{\Delta t}\left(\frac{4}{3}\boldsymbol{u}_g^{n-1} - \frac{1}{3}\boldsymbol{u}_g^{n-2}\right) + \boldsymbol{F}_g^n\right) \end{aligned} \tag{2.13}$$

We emphasize that the coefficients $w_{K,f,g}$ for the absorbing boundary conditions are calculated as for Robin boundaries with Robin weights $\boldsymbol{\mathcal{D}}_{\Delta t}$.

### 2.2.2 The Newmark Method

For discretization in time on the internal part of the domain we use the second order accurate Newmark method, which is widely used for solid mechanics. The original formulation of the Newmark method reads (Newmark, 1959):

$$\dot{u}^n = \dot{u}^{n-1} + (1-\gamma)\,\Delta t\ddot{u}^{n-1} + \gamma\Delta t\ddot{u}^n, \tag{2.14}$$

$$u^n = u^{n-1} + \Delta t\,\dot{u}^{n-1} + \frac{\Delta t^2}{2}[(1-2\beta)\ddot{u}^{n-1} + 2\beta\ddot{u}^n)], \tag{2.15}$$

where $\beta$ and $\gamma$ are the Newmark discretization parameters.

The primary variable of the MPSA discretization is displacement, meaning that we solve the linear system with displacements as the unknowns. The Newmark method can be adapted to suit the formulation of displacements as the only unknowns by rearranging Equations (2.14) and (2.15) (Chopra, 2012). This provides us with the following expressions for the current time-step velocity and acceleration:

$$\dot{u}^n = \left(1 - \frac{\gamma}{\beta}\right)\dot{u}^{n-1} + \Delta t\left[1 - \frac{\gamma}{2\beta}\right]\ddot{u}^{n-1} + \frac{\gamma}{\beta\Delta t}(u^n - u^{n-1}), \tag{2.16}$$

$$\ddot{u}^n = \frac{1}{\beta\Delta t^2}\left[u^n - u^{n-1} - \Delta t\dot{u}^{n-1} - (1-2\beta)\frac{\Delta t^2}{2}\ddot{u}^{n-1}\right]. \tag{2.17}$$

### 2.2.3 The Fully Discretized Model

The fully discretized model is obtained by combining the spatial and temporal discretizations presented in Sections 2.2.1 and 2.2.2. Substituting the discrete numerical tractions (Equation (2.13)) and the Newmark acceleration (Equation (2.17)) into the space-discrete system of ordinary differential equations represented by Equation (2.10) provides the fully discrete elastic wave equation and boundary conditions with displacements as the unknowns. To close the system of discretized equations, we assign initial boundary displacement values for the two previous time-steps (see Equation (2.13)), and initial values of displacement, velocity and acceleration on the internal part of the domain on the previous time-step. After initialization, the fully discretized model equations are solved for the displacements. After the displacements are obtained, the velocity and acceleration are updated according to Equations (2.16) and (2.17), respectively.

# 3 Analysis of the MPSA-Newmark Method

This section presents analysis of the MPSA-Newmark method with different types of boundary conditions. The first subsection is dedicated to convergence analysis of the MPSA-Newmark method with Dirichlet boundary conditions, including both theoretical considerations and a numerical convergence study. The second subsection contains a numerical convergence analysis and a numerical stability analysis of the MPSA-Newmark method with absorbing boundary conditions.

We now define the following $L^2$-norms for the cell quantity $\boldsymbol{u}$ and the face quantity $\boldsymbol{T}$:

$$\|\boldsymbol{u}\|_{\mathcal{T}} = \left(\sum_{K\in\mathcal{T}} m_K(\boldsymbol{u}_K \cdot \boldsymbol{u}_K)\right)^{\frac{1}{2}}, \tag{3.1}$$

$$\|\boldsymbol{T}\|_{\mathcal{F}} = \left(\sum_{f\in\mathcal{F}} \frac{1}{D} m_f(d_L - d_R)\cdot n_f\left(\boldsymbol{T}_f \cdot \boldsymbol{T}_f\right)^2\right)^{\frac{1}{2}}, \tag{3.2}$$

where $L, R \in \mathcal{T}_f$ in such a way that $\mathcal{F}_L \cap \mathcal{F}_R \neq \emptyset$. The symbol $d_L$ denotes the vector pointing from the face-center $x_f$ of $f$ to the cell-center $x_L$ of $L$. Similarly, $d_R$ is the vector pointing from $x_f$ to $x_R$. The symbol $n_f$ denotes the face normal vector pointing from cell $R$ to cell $L$, and the symbol $D$ represents the dimension of the simulation domain. In the case of $f \in \tilde{\mathcal{F}}$, one of $d_L$, $d_R$ is the zero vector.

The numerical errors of cell- and face-centered quantities are computed relative to the exact solutions projected onto cell-centers and face-centers. Hence, we introduce the projection $\Pi_{\mathcal{T}}$ which returns a cell-centered quantity $\left((\Pi_{\mathcal{T}}u)_K = u(x_K)\right)$ and the projection $\Pi_{\mathcal{F}}$ which returns a face-centered quantity $\left((\Pi_{\mathcal{F}}T)_f = T(x_f)\right)$. Then, considering plain type symbols to denote exact analytical solutions and boldface symbols to denote discrete solutions, the relative $L^2$-errors $\mathcal{E}$ are computed according to:

$$\mathcal{E}_u = \frac{\|\boldsymbol{u} - \Pi_{\mathcal{T}}u\|_{\mathcal{T}}}{\|\Pi_{\mathcal{T}}u\|_{\mathcal{T}}}, \tag{3.3}$$

$$\mathcal{E}_T = \frac{\|\boldsymbol{T} - \Pi_{\mathcal{F}}T\|_{\mathcal{F}}}{\|\Pi_{\mathcal{F}}T\|_{\mathcal{F}}}. \tag{3.4}$$

Subscript $u$ and $T$ on the error symbol denote displacement and traction errors, respectively. All errors are computed at the final time.

The linear systems in this section that arise from the 2D simulations are solved by a direct solver. The linear systems from the 3D simulations are solved by GMRES with the Geometric Algebraic Multigrid (GAMG) preconditioner in *petsc4py* (Dalcin et al., 2011). Runscripts for generating the numerical results presented in the following sections are accessible in the Docker image found in Jacobsen et al. (2025).

## 3.1 Convergence of MPSA-Newmark with Dirichlet Boundary Conditions

The combination of MPSA and Newmark can be understood in the context of the so-called "Method of Lines" (MoL). Thus, we consider the problem after spatial discretization as a semi-discrete problem (discrete in space and continuous in time), to which we thereafter apply the temporal discretization to resolve the time evolution. This perspective is advantageous, as the MoL is known to converge to the full solution of the time-space problem under quite weak conditions for the spatial and temporal discretization (see e.g. Knabner and Angermann (2003), Verwer and Sanz-Serna (1984)). Relevant for the current study is that for the space-time discretization to be convergent, the temporal discretization needs to be C-stable (see Definition 4.1 in (Verwer & Sanz-Serna, 1984)) independent of the stiffness of the problem arising from the spatial discretization (Verwer & Sanz-Serna, 1984). A common choice for stable discretization of second-order derivatives in time is the Newmark method, and in particular with the parameter choices $\beta = 1/4$ and $\gamma = 1/2$. This parameter choice provides a stable and energy conserving scheme (see e.g. (Subbaraj & Dokainish, 1989)), and is also the choice we employ in the current work. Due to the linearity of our model problem, the convergence rate will be limited by the convergence rate of the spatial and temporal discretization schemes. In particular, let $\psi$ denote convergence rate of the temporal method. The spatial convergence of the primary variable, namely the displacement $u$, is denoted by $\phi_u$, while that of the secondary variable, namely the traction $T$, is denoted by $\phi_T$. We then expect the space-time approximation, where all errors are computed at the final time, to satisfy:

$$\begin{aligned} \|\boldsymbol{u} - \Pi_{\mathcal{T}}u\|_{\mathcal{T}} &= \mathcal{O}\left((\Delta x)^{\phi_u} + (\Delta t)^{\psi}\right), \\ \|\boldsymbol{T} - \Pi_{\mathcal{F}}T\|_{\mathcal{F}} &= \mathcal{O}\left((\Delta x)^{\phi_T} + (\Delta t)^{\psi}\right), \end{aligned} \tag{3.5}$$

where $\Delta x$ denotes the grid size. For the MPSA method, the convergence rates are established numerically to be $\phi_u = 2$ and $1 < \phi_T \leq 2$ for a broad class of grids (Keilegavlen &

Nordbotten, 2017; Nordbotten & Keilegavlen, 2021). For the Newmark method with the abovementioned parameter choice, the convergence rate is $\psi = 2$ (Subbaraj & Dokainish, 1989). The temporal convergence rate of $\psi = 2$ is expected to hold for both displacements and tractions. We note that the spatial traction error will dominate in the case of refining with equal rates in both space and time, resulting in a lower convergence of the tractions than that of the displacements in a combined space-time convergence analysis.

Numerical convergence analysis of the MPSA-Newmark with Dirichlet boundary conditions is performed in three dimensions for unstructured simplex grids based on a known analytical reference solution. We performed the analysis on a unit cube domain, and all material parameters, including the density, hold unitary values corresponding to an isotropic solid. The body force term, $q$, in Equation (2.1) as well as initial values for displacement, velocity and acceleration are chosen in accordance with the analytical solution. The analytical solution holds zero Dirichlet boundary conditions. The analytical solution of the displacement for the convergence analysis is $u = \left[u_x, u_y, u_z\right]^T$, where each component of the solution $u$ is chosen as follows:

$$
\begin{aligned}
u_x &= \sin\left(\frac{5\pi}{2}t\right)\sin(\pi x)\sin(\pi y)\sin(\pi z) \\
u_y &= \sin\left(\frac{5\pi}{2}t\right)xyz(1-x)(1-y)(1-z) \\
u_z &= \sin\left(\frac{5\pi}{2}t\right)xy(1-x)(1-y)\sin(\pi z)
\end{aligned}
\tag{3.6}
$$

An illustration of the discrete solution at the final time ($t = 1.0$ s) is presented in Figure 3.1.

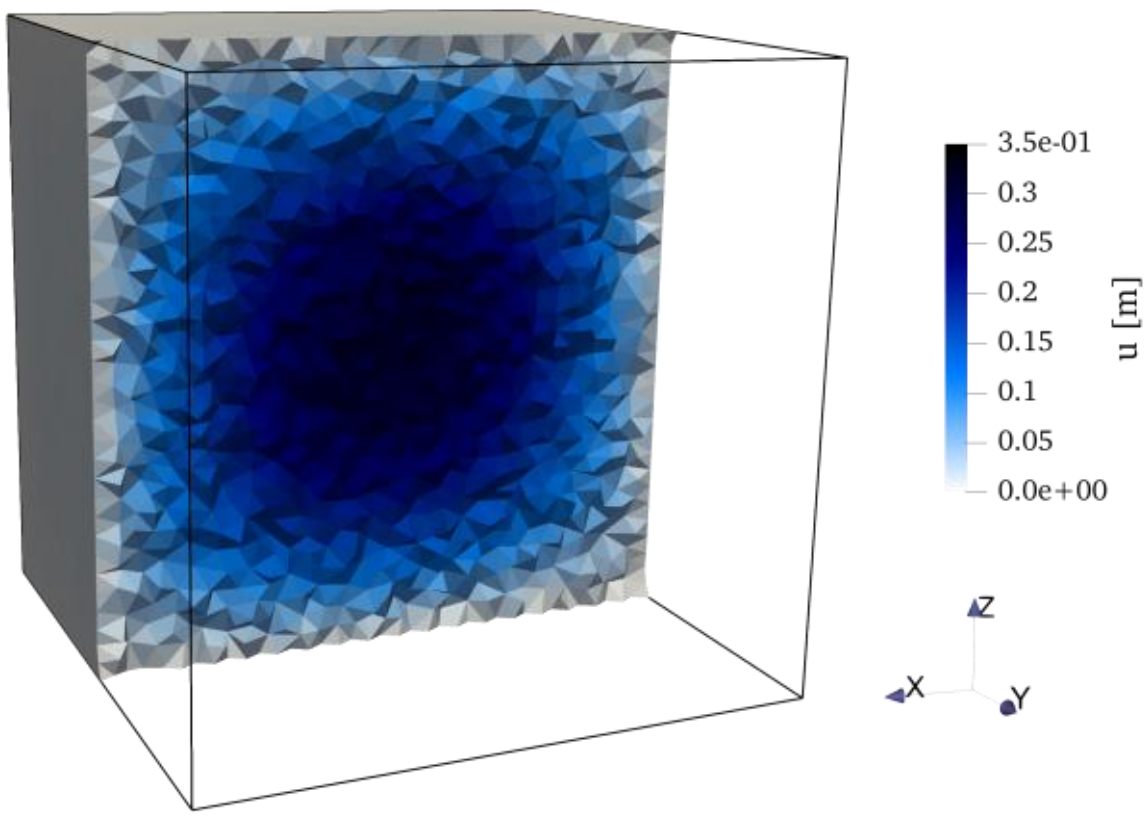

*Figure 3.1: Illustration of the displacement magnitude at time $t = 1.0\ s$ for the finest grid in the convergence analysis. The domain is cut along $y = 0.5\ m$ to visualize a cross-section of the displacement profile.*

The convergence analysis is performed with successive refinement by a factor of 2 in both space and time. The results of the analysis are displayed in Table 3.1 and Figure 3.2 with respect to number of cells $N_x$ and number of time-steps $N_t$. We see that the displacement error converges with 2nd order while the traction error converges with between 1st and 2nd order, which are the expected rates. We observe the same convergence rates also for Cartesian grids, but for the sake of brevity, the error values and plot are omitted from the text. The runscript which can reproduces results for Cartesian grids, as well as the results for the simplex grid, is found in the Docker Image in Jacobsen et al. (2025).

*Table 3.1: Table for refinement parameters and error values for the combined space-time convergence analysis. The table presents the number of cells, number of time-steps and error of displacement and traction for all refinement stages.*

| Refinement | $N_x$ | $N_t$ | Displacement error | Traction error |
|---|---|---|---|---|
| Initial | 387 | 150 | 3.01e-01 | 2.79e-01 |
| First | 2570 | 300 | 4.57e-02 | 6.55e-02 |
| Second | 19071 | 600 | 1.06e-02 | 1.99e-02 |
| Final | 148129 | 1200 | 2.66e-03 | 7.16e-03 |

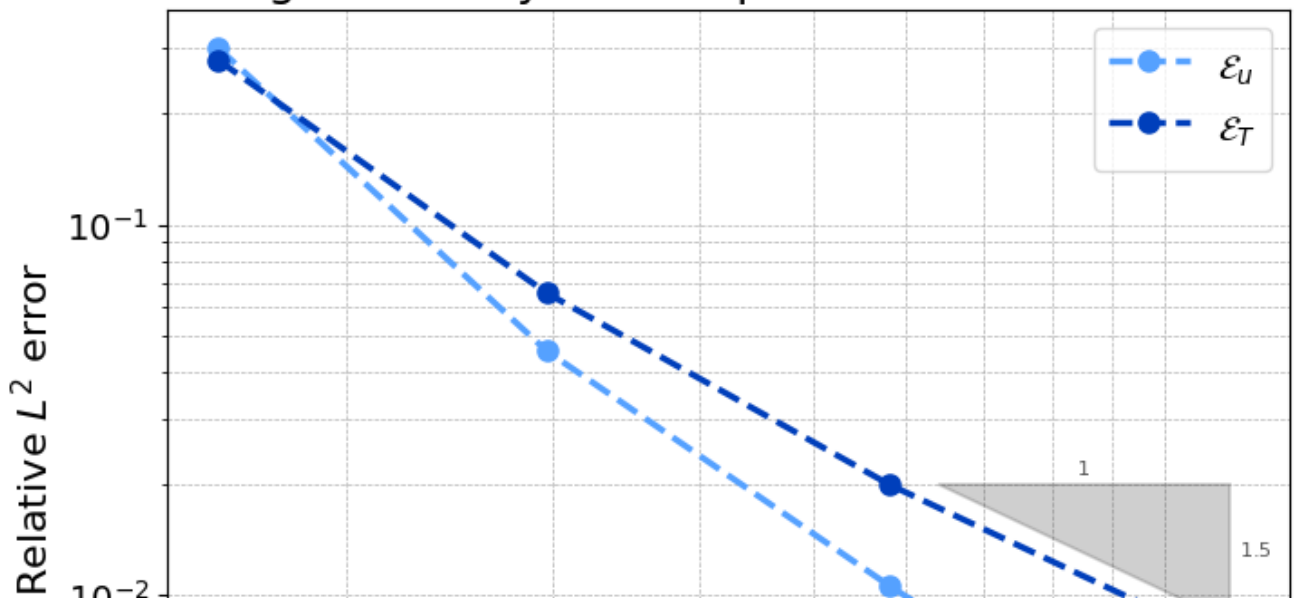

*Figure 3.2: Error plot for the space-time convergence analysis of MPSA-Newmark with Dirichlet boundaries. The $x$-axis is a combined space-time parameter consisting of the number of cells $N_x$ and the number of time-steps $N_t$.*

## 3.2 MPSA-Newmark with Absorbing Boundary Conditions

The absorbing boundary conditions introduce a first-order derivative in time which we discretize by a second-order implicit method. Thus, the temporal discretization of the problem becomes a combination of Newmark and the second-order implicit method. Therefore, the convergence rates for Dirichlet boundary conditions cannot necessarily be expected to carry over to a model with absorbing boundary conditions. In this section we therefore present numerical convergence analyses of the MPSA-Newmark scheme for a setup with absorbing boundary conditions. Additionally, to demonstrate the stability of the scheme with the inclusion of absorbing boundaries, we present an energy decay analysis.

### 3.2.1 Convergence Analysis

As outlined in the introduction, numerical methods which can accurately solve problems in heterogeneous and anisotropic media are highly relevant for applications such as extraction of geothermal energy, $CO_2$-storage and wastewater disposal. Therefore, we present a numerical convergence analysis of the MPSA-Newmark in homogeneous, heterogeneous, isotropic and anisotropic media, as well as combinations thereof.

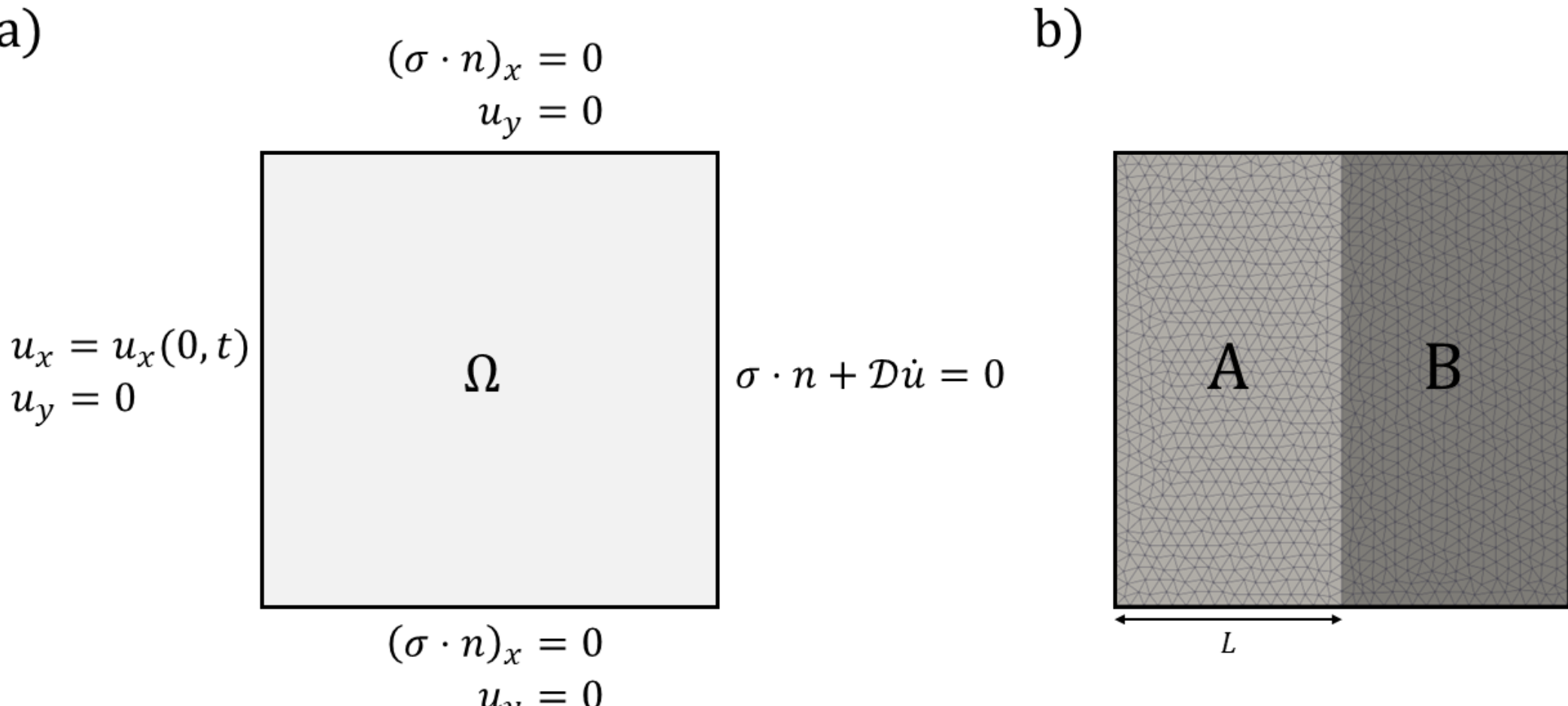

*Figure 3.3: a) Schematic of the boundary conditions for the convergence analyses of MPSA-Newmark with inclusion of absorbing boundary conditions. Subscript $x$ and $y$ on quantities denote the $x$ and $y$ component of the quantity. b) Schematic of the domain used in the convergence analyses. The domain is split vertically at $x = L$. Domain A is homogeneous isotropic, while domain B is either homogeneous and isotropic, or homogeneous and anisotropic in the $y$-direction.*

The analysis is performed on 2D unit squares with a known solution $u(x, t)$ which represents a wave that travels with normal incidence angle towards the right. For all the setups we have a time-dependent Dirichlet condition on the left side, an absorbing condition on the right side and rolling conditions on the top and bottom. We refer to Figure 3.3 a) for a schematic of the domain and its boundary conditions.

We refer to Figure 3.3 b) for a schematic of how $\Omega$ is partitioned into two regions A and B, as well as an illustration of how the domain is discretized using an unstructured simplex grid which conforms to the boundary between A and B. The different configurations of A and B is what determines whether the setup is homogeneous, heterogeneous, isotropic or anisotropic.

The 2D stiffness tensors used in this section can be derived from Equation (2.5). This is done by choosing $v$ to align with the $y$-axis, as we are only considering anisotropy in the direction orthogonal to the wave. Additionally, we set $\lambda_{||} = 0$ and $\mu = \mu_{||} = \mu_{\perp}$. The stiffness tensors are consequently determined by three independent material parameters in the following way:

$$C_{ijkl} = \lambda \delta_{ij}\delta_{kl} + \mu\left(\delta_{ik}\delta_{jl} + \delta_{il}\delta_{jk}\right) + \lambda_{\perp}\delta_{i2}\delta_{j2}\delta_{k2}\delta_{l2}, \quad i, j, k, l = 1, 2. \tag{3.7}$$

Here $\lambda_{\perp}$ is what determines whether Equation (3.7) represents an isotropic or an anisotropic medium. In domain A, which is always isotropic, we set $\lambda_{\perp} = 0$. In domain B, which is either isotropic or anisotropic, we set $\lambda_{\perp} \geq 0$.

Common for all the setups is that we have the following analytical solution for the wave in the medium:

$$u(x,t) = \begin{cases} \left[\sin\left(t - \dfrac{x - L}{c_A}\right) + \dfrac{c_A - c_B}{c_A + c_B}\sin\left(t + \dfrac{x - L}{c_A}\right), 0\right]^{\mathrm{T}} & \text{for } x < L, \\ \left[\dfrac{2c_A}{c_B + c_A}\sin\left(t - \dfrac{x - L}{c_B}\right), 0\right]^{\mathrm{T}} & \text{for } x \geq L. \end{cases} \tag{3.8}$$

Here $L$ denotes the width of domain A, and thus the location (on the $x$-axis) of the boundary between A and B. $c_A$ and $c_B$ represent the primary wave speed in A and B, respectively, where the expression for the primary wave speed is given by $c = \sqrt{\frac{\lambda+2\mu}{\rho}}$. For all the simulation setups we set initial displacement values according to Equation (3.8), and initial velocity and acceleration values according to the first and second time derivative of Equation (3.8), respectively.

The different configurations of domains A and B are determined by the value of the material parameters in Equation (3.7), as well as the heterogeneity coefficient, $r_h$, and the anisotropy coefficient, $r_a$. The heterogeneity and anisotropy coefficients determine whether domain B is like A, if it is isotropic with different parameter values than A or if it is anisotropic. We refer to Table 3.2 for an overview of the parameter values we have employed. Values of the heterogeneity and anisotropy coefficients are presented alongside the convergence results for the respective convergence runs, and values for space and time related quantities are found in Table 3.3.

*Table 3.2: Parameter values for domain* A *and* B.

| Symbol | Values for domain A | Values for domain B | Unit |
|---|---|---|---|
| $\lambda$ | 0.01 | $0.01r_h$ | $\text{kg m}^{-1}\text{s}^{-2}$ |
| $\mu$ | 0.01 | $0.01r_h$ | $\text{kg m}^{-1}\text{s}^{-2}$ |
| $\lambda_\perp$ | 0.0 | $0.01r_a$ | $\text{kg m}^{-1}\text{s}^{-2}$ |
| $\rho$ | 1.0 | 1.0 | $\text{kg m}^{-3}$ |

In Figure 3.4 we present convergence analysis results of MPSA-Newmark with absorbing boundary conditions for 6 different configurations of A and B. The method converges as expected for isotropic, homogeneous and heterogeneous media, as well as media where there is anisotropy in the direction orthogonal to the propagating wave. The figure shows that the method converges for both low and high heterogeneity coefficients, but that higher heterogeneity coefficients result in larger errors. The figure also shows that errors from the anisotropic setups ($r_a = 10$) mostly overlap with those from the isotropic setups ($r_a = 0$) for all three different heterogeneity ratios. This indicates that the numerical error in the wave propagation is not influenced by the anisotropy, even for large parameter contrasts ($\frac{r_a}{r_h} \approx 2500$), showing the method's robustness in handling anisotropy.

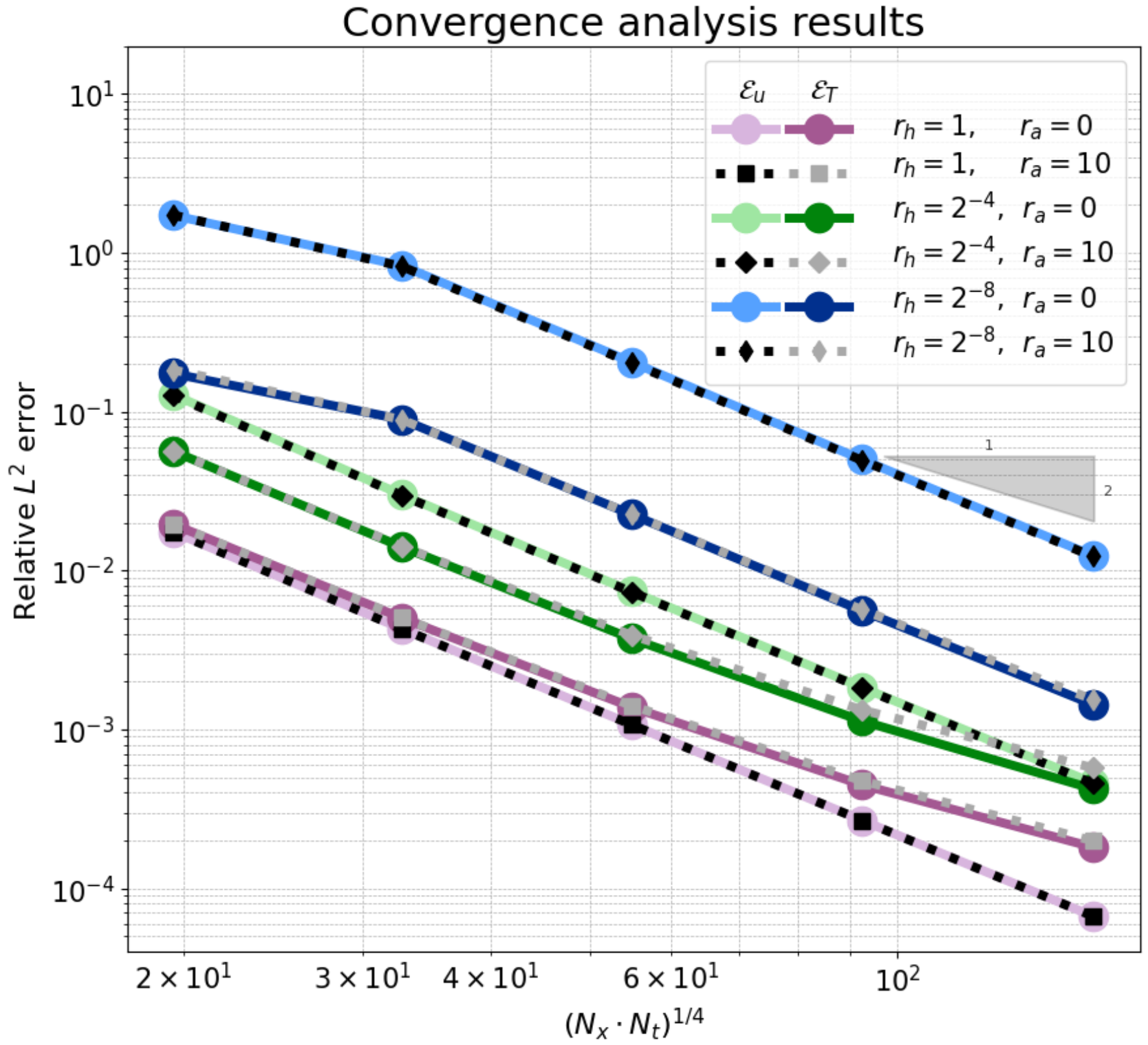

*Figure 3.4: Displacement and traction errors from the combined temporal and spatial convergence analysis of MPSA-Newmark with absorbing boundaries. The errors obtained for the setups with vertical anisotropy have quadrilateral markers and are colored in black and grey. The errors obtained for the isotropic setups have circular markers and are given colors. The $x$-axis of the plot is represented by a combined space-time parameter consisting of the number of cells $N_x$ and the number of time-steps $N_t$.*

*Table 3.3: Time and space related parameters for the temporal and spatial convergence analysis of MPSA-Newmark in isotropic, heterogeneous and anisotropic media with absorbing boundary conditions.*

| | Value | Unit |
|---|---|---|
| $\Omega$ | $[0, 1.0]^2$ | $\mathrm{m}^2$ |
| $L$ | 0.5 | m |
| $t_{final}$ | 15.0 | s |
| Initial $N_x$ | 2442 | – |
| $N_x$, first refinement | 9568 | – |
| $N_x$, second refinement | 38014 | – |
| $N_x$, third refinement | 151802 | – |
| Final $N_x$ | 606672 | – |
| Initial $\Delta t$ | 0.25 | s |
| $\Delta t$, first refinement | 0.125 | s |
| $\Delta t$, second refinement | 0.0625 | s |
| $\Delta t$, third refinement | 0.03125 | s |
| Final $\Delta t$ | 0.015625 | s |

### 3.2.2 Energy Decay Analysis

We here demonstrate the stability of the scheme by considering propagation of waves with non-orthogonal incidence angles onto absorbing boundaries in an isotropic domain. See Table 3.4 at the end of this section for material and domain parameters. All the boundaries in the domain are absorbing, and as there are no driving forces present in the system, the energy within the domain is expected to decrease towards zero as time passes. The discrete energy is described by $E = \sum_{K\in\mathcal{T}} \int_K \rho(\dot{\boldsymbol{u}}_K \cdot \dot{\boldsymbol{u}}_K)\, dV$, where $\dot{\boldsymbol{u}}_K$ is the numerical velocity in cell $K$.

The wave is imposed by initial values for displacement, velocity and acceleration, where the expression for the initial displacement wave is as follows:

$$u(x,y,t=0) = \begin{bmatrix} \cos(\theta)\sin\left(-\dfrac{\cos(\theta)\,x + \sin(\theta)\,y}{c}\right) \\ \sin(\theta)\sin\left(-\dfrac{\cos(\theta)\,x + \sin(\theta)\,y}{c}\right) \end{bmatrix}, \tag{3.9}$$

We refer to Figure 3.5 for an illustration of the initial displacement profile for both an orthogonal wave and a rotation of the orthogonal wave:

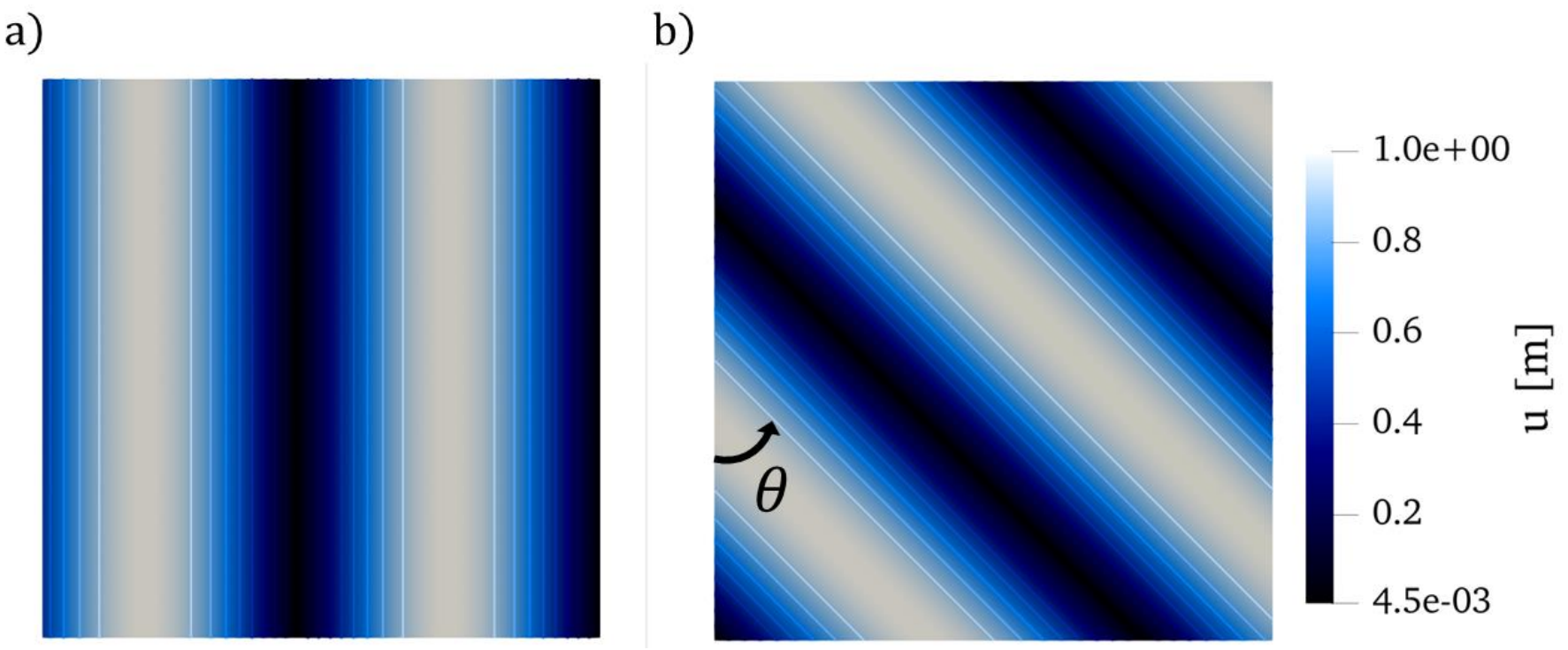

*Figure 3.5: Displacement magnitude for a) an orthogonal wave. b) a rotation of the orthogonal wave. Contour lines for the displacement magnitude are included for illustration purposes.*

The energy decay analysis is first performed for waves rotated by $\frac{\pi}{4}$ rad for various grid sizes. The simulation is run until a final time of $t = 15.0$ s with a time-step $\Delta t = 1/20$ s. Grid sizes and the corresponding number of cells are listed in Table 3.5 at the end of this section.

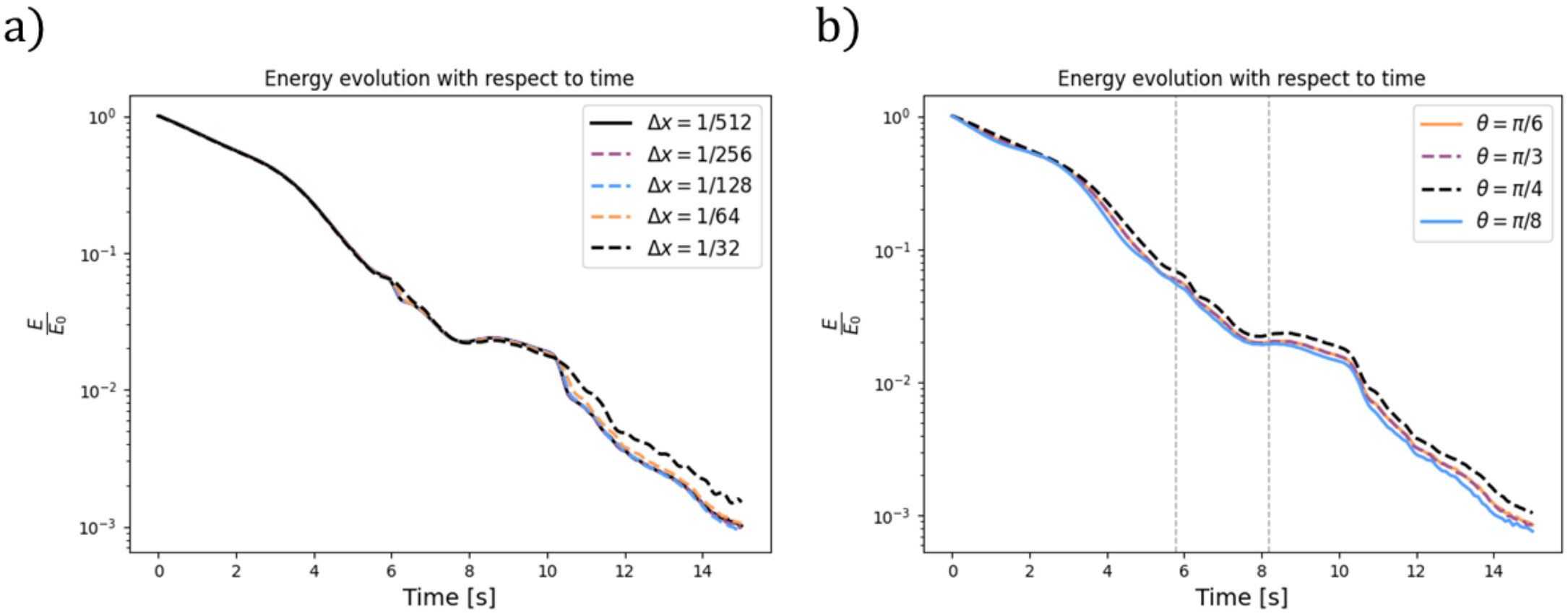

*Figure 3.6: a) Energy evolution as time passes for different grid sizes, $\Delta x$. The wave rotation angle is held constant at $\pi/4$. b) Energy evolution as time passes for different rotation angles $\theta$. The left dashed line corresponds to the time when an orthogonal wave has exited the domain, while the right one corresponds to when a $\frac{\pi}{4}$ rad rotated wave has exited the domain.*

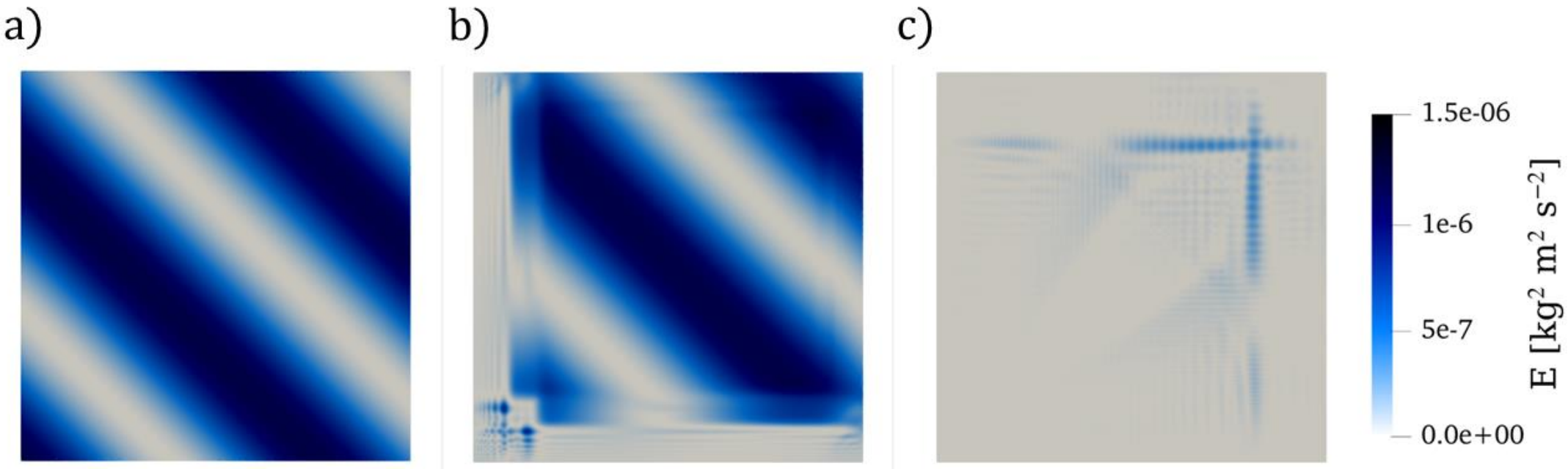

*Figure 3.7: Discrete energy distribution in a domain with $N_x = 606496$ cells. a) At initialization ($t = 0.0$ s), b) At $t = 1.0$ s. c) At $t = 8.5$ s, which is when the initial wave has exited and only the initialization artifact remains.*

In Figure 3.6 a) we show the temporal evolution of the ratio between the energy $E$ and the initial energy $E_0$ for successive grid refinement. The figure shows that the energy decay can be divided into three phases: Decreasing, then plateauing and finally decreasing again. Looking to Figure 3.7 a) we see the initial energy distribution for a rotated wave. As the wave propagates towards the upper right corner, no new wave is introduced to the domain, which in turn causes a discontinuity in the solution along the bottom and left boundary. This leads to the emergence of a spurious wave in the lower left corner, which is common when using higher order time-discretization schemes in the presence of discontinuities (see Figure 3.7 a) and b)). Referring to Figure 3.7 c), which is a snapshot from after the initialized wave has exited the domain, we can see that the spurious wave is still present due to its significantly lower propagation speed. The energy of the spurious wave is represented in the energy plateau between 8.0 s and 10.0 s in Figure 3.6. When the slow wave finally exits the domain after around 10.0 seconds, we see a rapid energy decrease.

Looking to Figure 3.6 b), we present the energy evolution of waves rotated by different angles $\theta$. The right dashed line corresponds to when a wave rotated by $\frac{\pi}{4}$ rad has exited the domain, and the left dashed line corresponds to when an orthogonal wave has exited domain. The region between the two dashed lines is thus the time interval where waves rotated at angles between 0 and $\frac{\pi}{4}$ rad are expected to have exited the domain. Values for the time-step size and number of cells for the simulations with different rotation angles correspond to the values labeled as "First refinement" in Table 3.4. The figure shows that the system energy decreases for all the values of $\theta$ that we tested.

The difference in energy decay for the different wave rotation angles is small, even though the waves should exit the domain at different times. This is in part due to the emergence of the spurious wave at the bottom left corner, but it is also likely due to the imperfectness of the absorbing boundary conditions we have employed.

*Table 3.4: Material parameters for the energy decay analysis.*

| Symbol | Value | Unit |
|---|---|---|
| $\lambda$ | 0.01 | kg m$^{-1}$ s$^{-2}$ |
| $\mu$ | 0.01 | kg m$^{-1}$ s$^{-2}$ |
| $\rho$ | 1.0 | kg m$^{-3}$ |
| $t_{\text{final}}$ | 15.0 | s |
| $\Omega$ | $[0, 1.0]^2$ | m |

*Table 3.5: Grid size, number of cells and time-step size for the energy decay analysis with successive grid refinement.*

| | Grid size, $\Delta x$ [m] | $N_x$ | $\Delta t$ [s] |
|---|---|---|---|
| Initial refinement | 1/32 | 2400 | 0.05 |
| First refinement | 1/64 | 9520 | 0.05 |
| Second refinement | 1/128 | 37976 | 0.05 |
| Third refinement | 1/256 | 151724 | 0.05 |
| Final refinement | 1/512 | 606 496 | 0.05 |

# 4 Wave Propagation: Simulation Examples

In this section we present simulation examples of wave propagation in anisotropic, heterogeneous and fractured media. The goal is to demonstrate that the methodology produces expected wave propagation behavior in different kinds of media. Expected behavior includes different wave propagation patterns in different directions in anisotropic media and reflections from boundaries between media with different stiffness. In the case of inclusion of a fracture, which in this section is modelled as open and not filled with any fluid, the expectation is that the waves are fully reflected.

All the linear systems for the simulation examples in this section are solved by GMRES with the GAMG preconditioner from *petsc4py* (Dalcin et al., 2011).

## 4.1 Example 1: Simulations with Different Seismic Source Locations

In this section we present wave propagation in anisotropic and heterogeneous media where the simulation domain $\Omega$ is a unit cube. The domain is heterogeneous in the sense that it contains both isotropic and anisotropic regions. The anisotropic region is located fully inside $\Omega$ and is denoted by $\omega$. We consider the anisotropic domain to be transversely isotropic with the $z$-axis as the axis of symmetry. The stiffness tensor for the anisotropic region is determined by five independent material parameters as presented in Equation (*2.5*) while the isotropic part of the domain has a stiffness tensor according to Equation (2.4).

The dimensions and location of the anisotropic region are specified for each simulation example. We refer to Figure 4.1 for an illustration of an example domain where the anisotropic region is located at the domain center and has side lengths of 0.5 m.

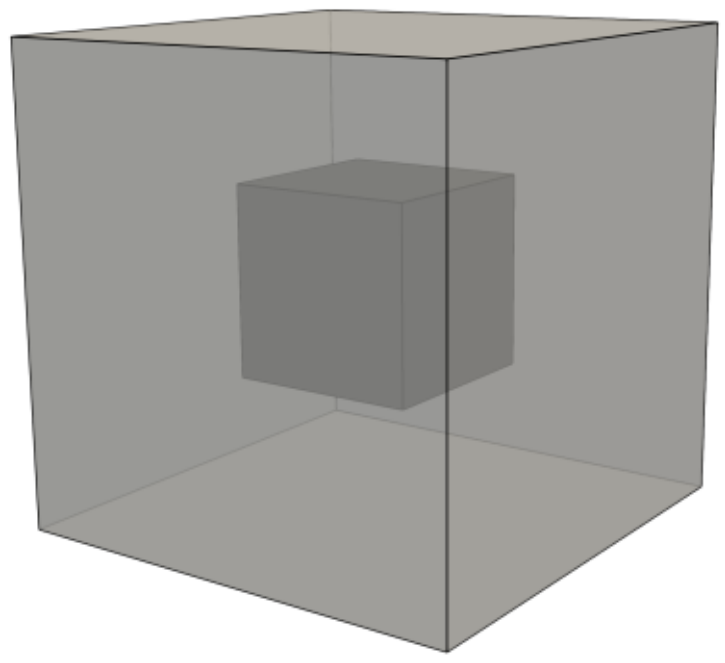

*Figure 4.1: Illustration of a 3D cubic domain with an inner domain. The lighter part of the domain is isotropic, while the darker, inner part is transversely isotropic. The location of the inner domain is different in the two simulation examples.*

Waves are introduced to the domain by a Ricker wavelet, which is commonly employed for modeling seismic sources. The simulation is initialized by zero initial displacement, zero initial acceleration and a velocity determined by the expression

$$\dot{u}(x, y, z) = e^{-\frac{\pi^2((x-x_R)^2+(y-y_R)^2+(z-z_R)^2)}{0.125^2}}[x - x_R, y - y_R, z - z_R]^T. \quad (4.1)$$

The velocity in Equation (4.1) corresponds to a Ricker wavelet, where $(\cdot)_R$ denotes the coordinate values for the center of the wavelet.

The dimensions of the outer and inner domain, as well as the material parameters for each of the domains, are held the same for both the simulation examples. The only difference between the two simulations is the location of the internal domain and the location of the wave source. An overview of the parameters for the simulation examples is presented in Table 4.1.

*Table 4.1 Simulation parameters for the heterogeneous and anisotropic simulation examples 1.1 and 1.2.*

| | Symbol | Value | Unit |
|---|---|---|---|
| **Common parameters** | | | |
| | $\Omega$ | $[0, 1.0]^3$ | $\mathrm{m}^3$ |
| | $N_x$ | 512 000 | – |
| | $N_t$ | 90 | – |
| | $\Delta t$ | 1/600 | s |
| | $t_{final}$ | 0.15 | s |
| **Isotropic domain parameters** | | | |
| | $\mu$ | 1.0 | $\mathrm{kg\ m^{-1}\ s^{-2}}$ |
| | $\lambda$ | 1.0 | $\mathrm{kg\ m^{-1}\ s^{-2}}$ |
| **Anisotropic domain parameters** | | | |
| | $v$ | $[0, 0, 1]^T$ | – |
| | $\mu_{\|\|}$ | 1.0 | $\mathrm{kg\ m^{-1}\ s^{-2}}$ |
| | $\mu_{\perp}$ | 2.0 | $\mathrm{kg\ m^{-1}\ s^{-2}}$ |
| | $\lambda_{\|\|}, \lambda_{\perp}$ | 5.0 | $\mathrm{kg\ m^{-1}\ s^{-2}}$ |
| | $\lambda$ | 1.0 | $\mathrm{kg\ m^{-1}\ s^{-2}}$ |
| **Location parameters for source within inner domain (Ex. 1.1)** | | | |
| | $\omega$ | $[0.25, 0.75]^3$ | $\mathrm{m}^3$ |
| | $(x_R, y_R, z_R)$ | $(0.5, 0.5, 0.5)$ | (m, m, m) |
| **Location parameters for source within outer domain (Ex. 1.2)** | | | |
| | $\omega$ | $[0.25, 0.75] \times [0.25, 0.75] \times [0.05, 0.55]$ | $\mathrm{m}^3$ |
| | $(x_R, y_R, z_R)$ | $(0.5, 0.5, 0.7)$ | (m, m, m) |

### 4.1.1 Example 1.1: Ricker source located within the anisotropic domain

We expect different wave propagation behaviors in isotropic and anisotropic domains. In isotropic materials we expect the waves to propagate uniformly in all directions, just like an expanding sphere. In anisotropic media the waves are expected to deviate from the radial symmetry that occurs in isotropic media, meaning that the wave can for instance have an ellipsoidal shape instead of a spherical one.

Using a transversely isotropic medium for the domain should allow us to observe both expected behaviors described above. The transversely isotropic medium considered herein has the $z$-axis as the axis of symmetry, meaning that the $xy$-plane is the plane of isotropy. The waves are thus expected to behave in an isotropic manner in the $xy$-plane, while all other planes are expected to hold distorted wave patterns.

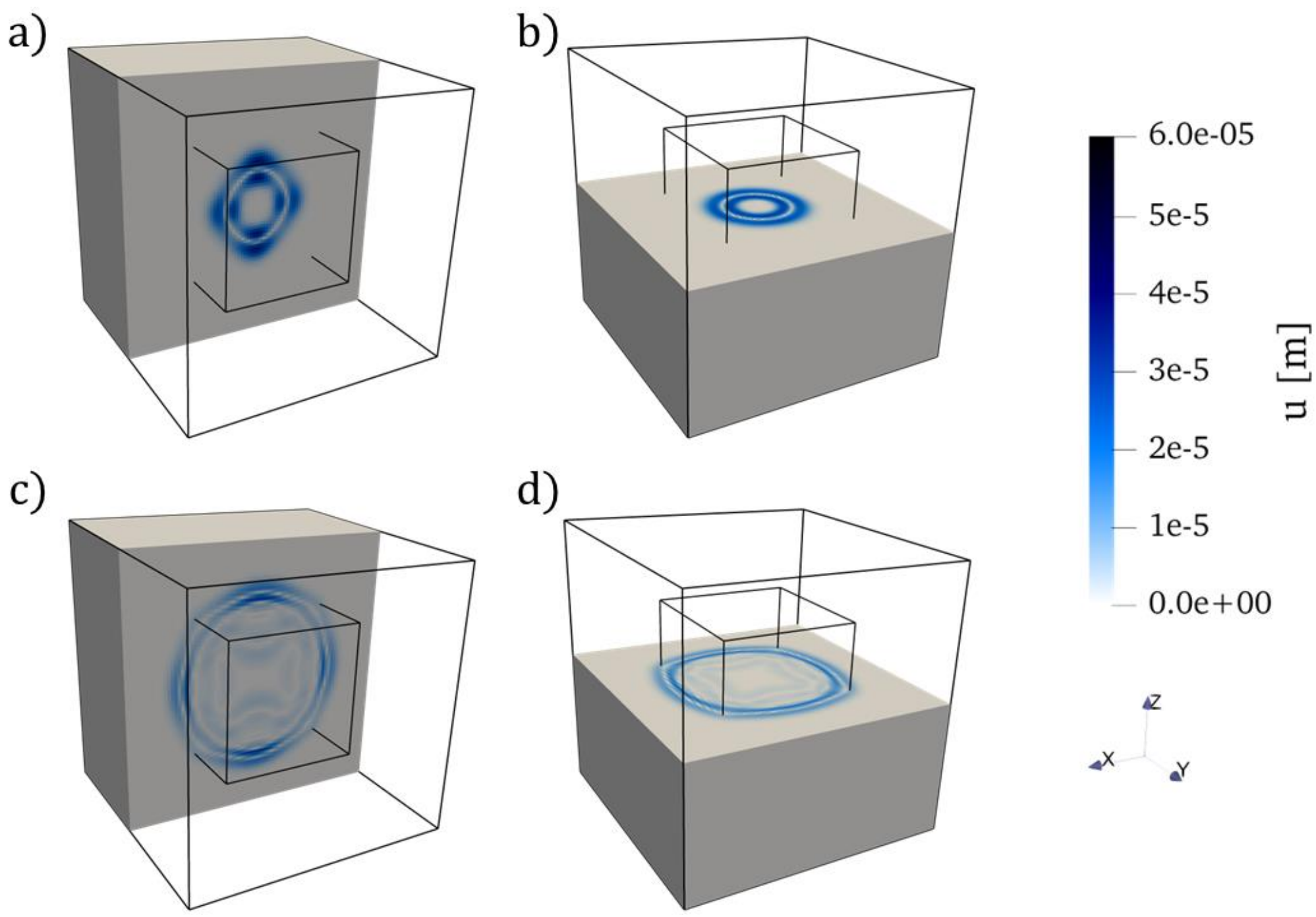

*Figure 4.2: Magnitude of the displacement vector for the simulation case where the propagating wave is initiated inside the anisotropic domain. The first row contains snapshots from time $t = 0.05$ s. The second row contains snapshots from time $t = 0.125$ s.*

Looking to Figure 4.2 a)-d) we see the wave propagation at two different times and in two different planes. The black outline within $\Omega$ represents the anisotropic region $\omega$. Figure 4.2 b) illustrates radial symmetry of the displacement wave, which is as expected as the figure depicts the displacement profile in the plane of isotropy at a point in time where the wave has not yet exited from the anisotropic region. On the other hand, Figure 4.2 a) and c) show only axial symmetry due to the anisotropies in the domain, while subfigure d) shows axial symmetry because the displacement wave now has propagated into the isotropic region. Additionally, subfigures c) and d) show that some reflections occur at the internal material discontinuity boundaries.

### 4.1.2 Example 1.2: Ricker source located outside the anisotropic domain

We will demonstrate elastic wave propagation in heterogeneous media further by locating the wave source outside the anisotropic region. As we saw in Section 4.1.1, there are wave reflections caused by intersections of material layers with different properties. When a wave travels towards and through heterogeneities, we also expect that both the shape and the propagation speed of the wave are altered. Depending on the values of the material parameters, the wave can move either faster or slower. In the case of the simulation presented in this section, we expect a higher propagation speed through the anisotropic domain due to higher values of $\lambda, \lambda_{||}, \lambda_{\perp}, \mu_{||}$ and $\mu_{\perp}$ (see Table 4.1).

As the Ricker wavelet is allowed to propagate through the isotropic medium, it eventually hits the anisotropic region which is indicated by the black outline within $\Omega$. In Figure 4.3 we see that the circular shape of the Ricker wavelet becomes distorted within the anisotropic region. Specifically, the wave travels quicker through the anisotropic region, just as expected from the choice of material parameters. In addition to this we see in Figure 4.3 c) and d) that a wave is reflected at the internal boundary.

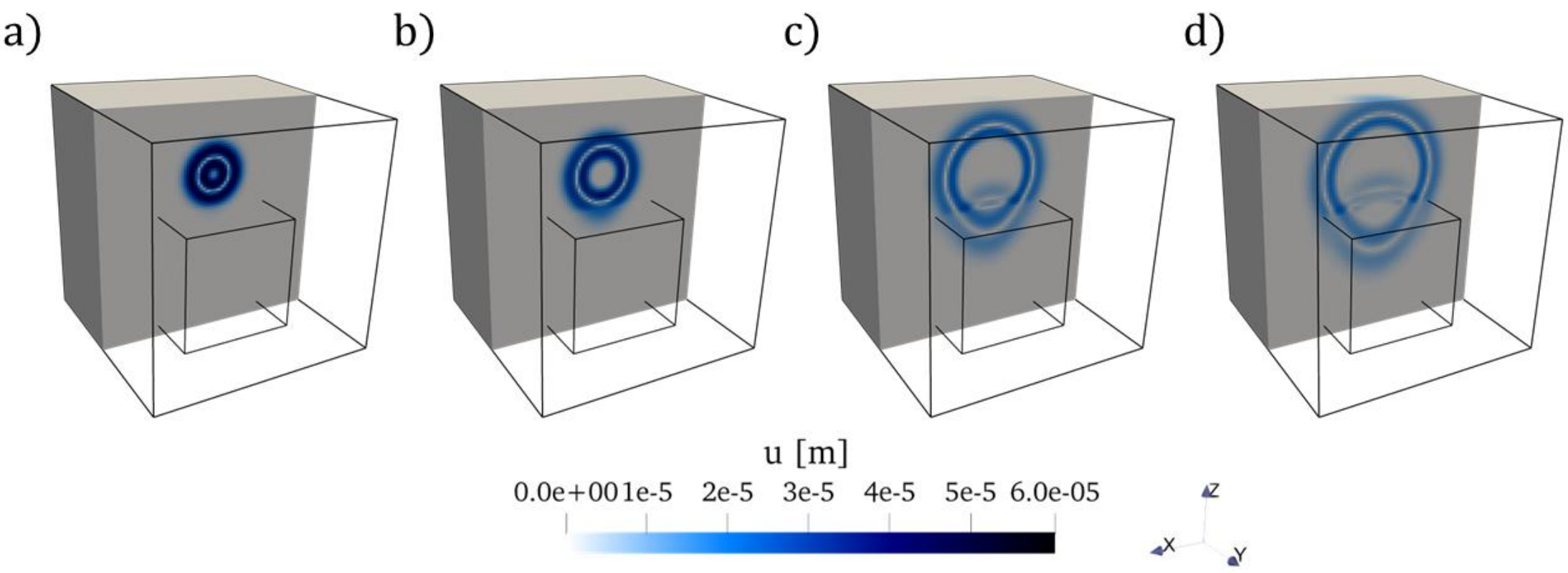

*Figure 4.3: Magnitude of the displacement vector for the simulation case where the propagating wave is initiated outside the anisotropic domain. The figures illustrate the Ricker wavelet's encounter with an internal domain boundary. Time a)* $t = 0.05$ s*, b)* $t = 0.075$ s*, c)* $t = 0.125$ s*, d)* $t = 0.15$ s*.*

## 4.2 Example 2: Simulation in Heterogeneous and Fractured Media

As outlined in the introduction, we are interested in applications related to the subsurface where the rock formations may be heterogeneous and fractured. Therefore, we present this example as a demonstration that we can use the present methodology, namely the MPSA-Newmark with absorbing boundary conditions, to investigate problems in such media.

The domain is divided into three isotropic layers with different material parameters. Parameters for the upper, middle and lower layer are presented in Table 4.2. We enforce absorbing boundary conditions on all the outer domain boundaries, and inside the domain we have an internal Neumann boundary with zero traction ($\sigma \cdot n = 0$) which models an open, empty fracture. The open fracture, which is a rectangular plane with corners in the points $(0.2, 0.2, 0.8)$, $(0.2, 0.8, 0.8)$, $(0.8, 0.8, 0.2)$ and $(0.8, 0.2, 0.2)$, is indicated in blue in the schematic of the domain in Figure 4.4.

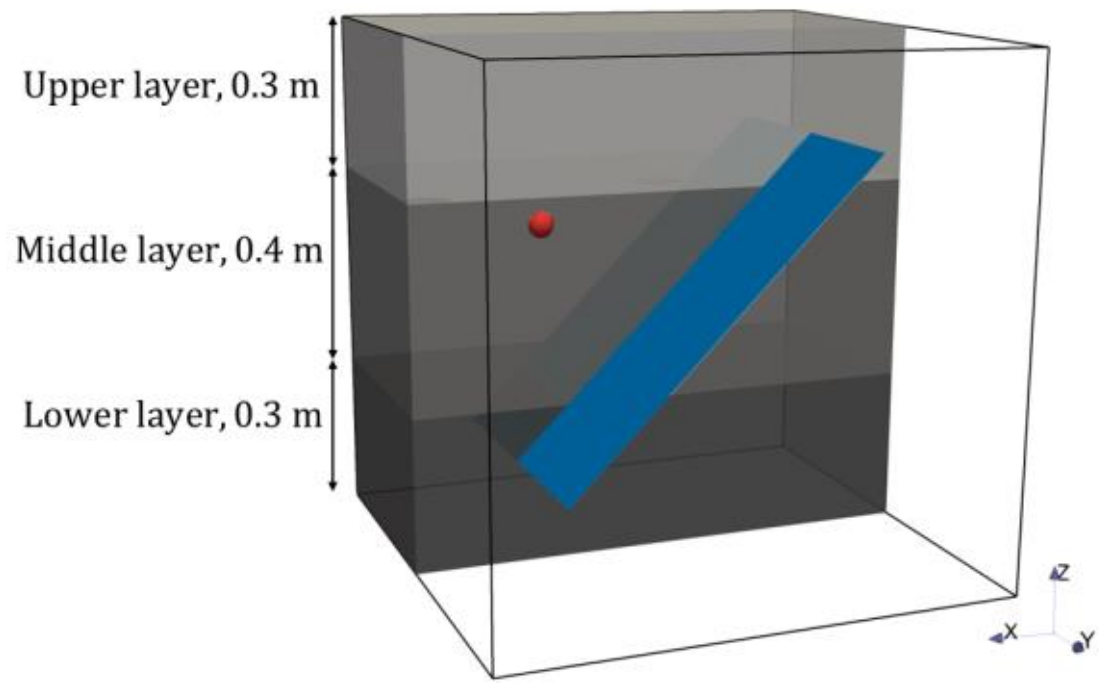

*Figure 4.4: Schematic of fractured and layered heterogeneous medium used in example 2. The red sphere indicates where the seismic source is located, and the blue plane represents an open fracture. The domain is cut at $y = 0.5$ m such that we get a cross-section view of the inside of the medium. Partial transparency is applied to show the fracture plane within the rock matrix.*

Like in Section 4.1 we stimulate the domain by a Ricker wavelet. That is, zero initial displacement and acceleration, and a velocity corresponding to the Ricker wavelet expression shown below.

$$\dot{u}(x, y, z) = e^{-\frac{\pi^2((x-0.75)^2+(y-0.5)^2+(z-0.65)^2)}{0.3^2}}[x - 0.75, y - 0.5, z - 0.65]^{\mathrm{T}}. \quad (4.2)$$

This corresponds to a Ricker wavelet centered at the point $(x, y, z) = (0.75, 0.5, 0.65)$. The location of the Ricker wavelet center is indicated by a red sphere in Figure 4.4.

We expect that the wave is fully reflected when it hits the open fracture, meaning that we expect seeing a sudden drop in displacement magnitude just beyond where the fracture is located. As the medium consists of layers of different parameter values, we also expect some reflections from the internal layer boundaries as well as different wave propagation speeds in the different layers, like that in Figure 4.3.

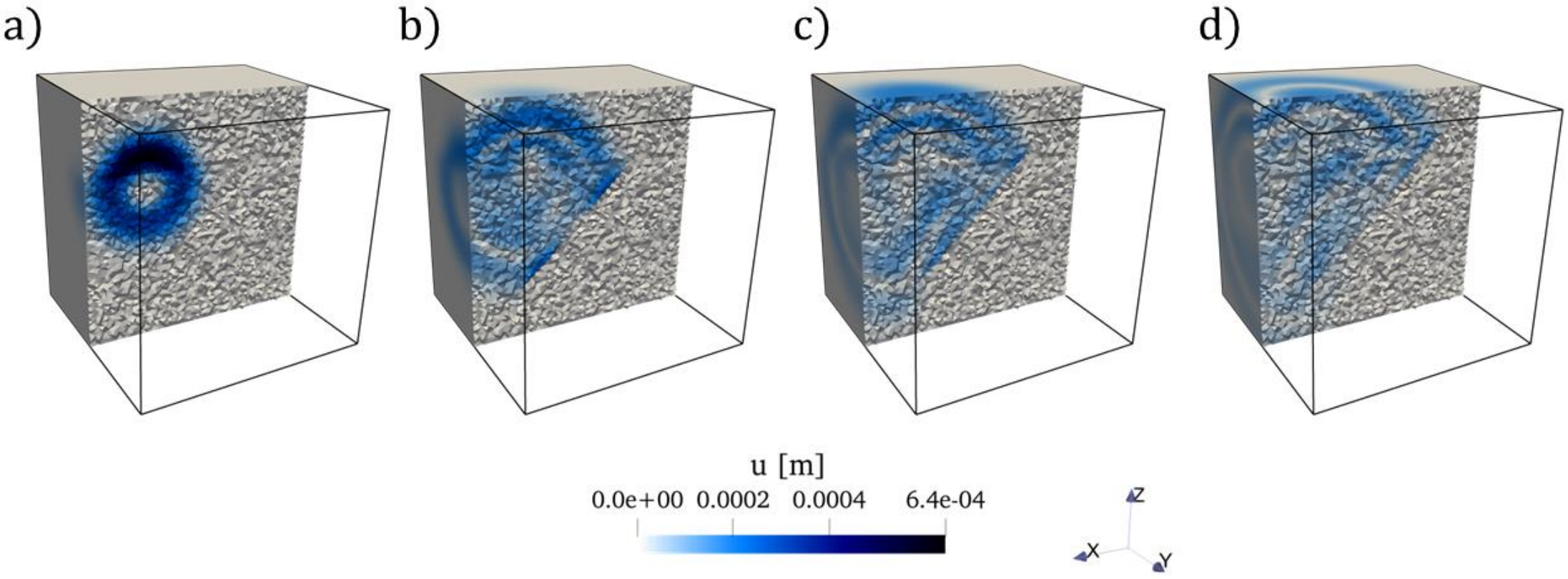

*Figure 4.5: Displacement vector magnitude for wave propagation in layered heterogeneous and fractured medium at four different times: a) $t = 0.05$ s, b) $t = 0.125$ s, c) $t = 0.175$ s, d) $t = 0.225$ s.*

Figure 4.5 illustrates how the wave propagates in the heterogeneous and fractured domain. Specifically, the figure shows that the wave is reflected as it hits the open fracture. Additionally, the outer domain boundaries allow for the wave to propagate outwards, leaving no visible artificial boundary reflections. Figure 4.5 also shows, especially in the $yz$-plane in subfigure c), wider waves in the lower layer. This is expected from the higher material parameters in this region than those found in the two other domain layers.

To illustrate the sudden drop in displacement magnitude across the fracture, we have plotted the displacement along the diagonal black line that is shown in Figure 4.6 b). The displacement along the line is presented in Figure 4.6 a), and it is shown that the displacement magnitude suddenly drops to zero at a distance of $1/\sqrt{2}$ m along the diagonal line. This is where the fracture is located, and the displacement drop indicates that the wave is fully reflected by the fracture.

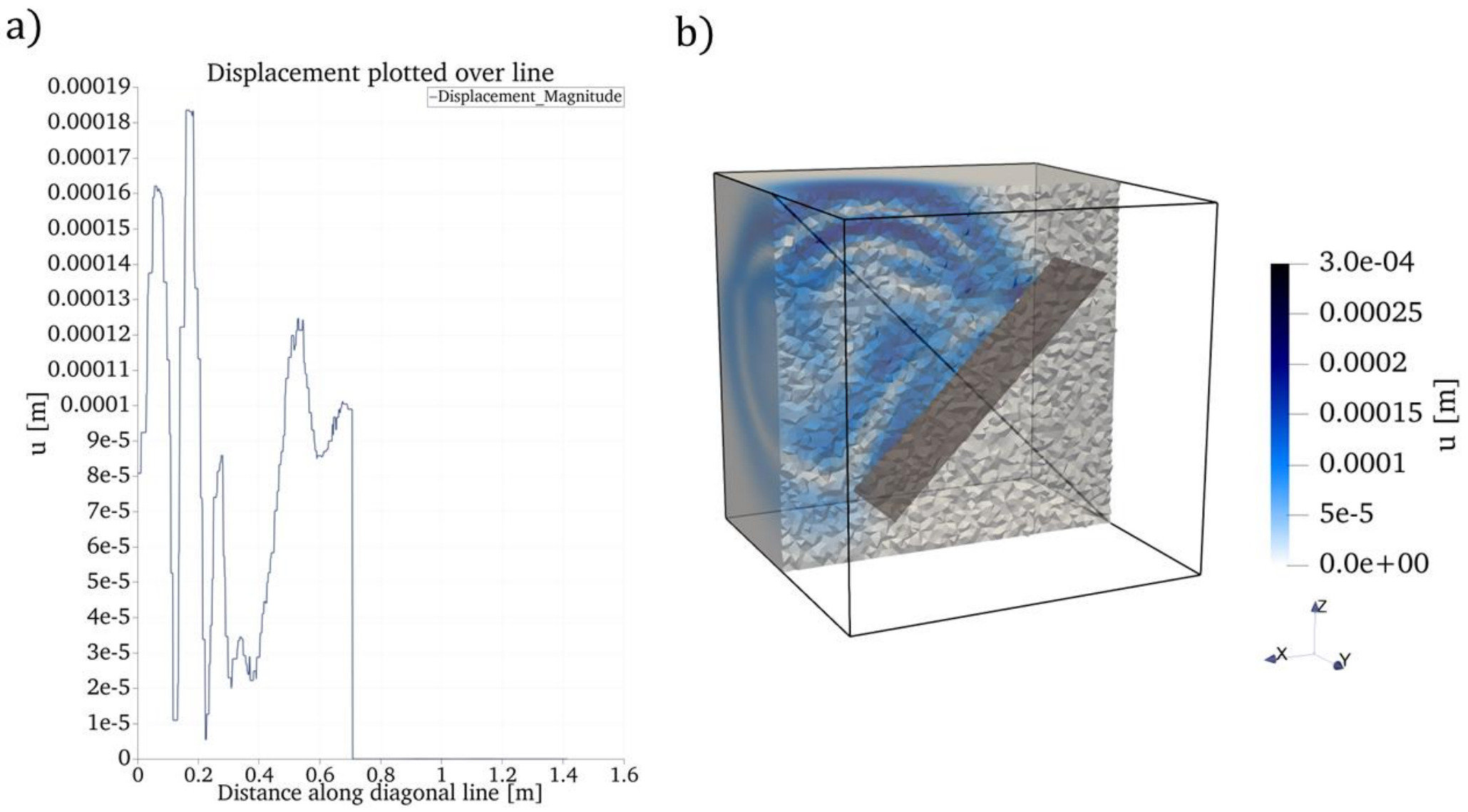

*Figure 4.6: Displacement magnitude for wave propagation in a fractured domain. The wave is reflected by the fracture at time $t = 0.175$ s. a) The displacement magnitude along the black diagonal line shown in subfigure b). b) Simulation domain cut at $y = 0.5$ m to visualize the wave propagation inside the domain. Partial transparency of the domain is applied to better visualize the dark grey fracture plane and the black diagonal line.*

*Table 4.2: Simulation parameters for example 2.*

| | Symbol | Value | Unit |
|---|---|---|---|
| **Upper layer** | | | |
| | $\lambda$ | 1.0 | $\mathrm{kg\ m^{-1}\ s^{-2}}$ |
| | $\mu$ | 1.0 | $\mathrm{kg\ m^{-1}\ s^{-2}}$ |
| **Middle layer** | | | |
| | $\lambda$ | 2.0 | $\mathrm{kg\ m^{-1}\ s^{-2}}$ |
| | $\mu$ | 2.0 | $\mathrm{kg\ m^{-1}\ s^{-2}}$ |
| **Lower layer** | | | |
| | $\lambda$ | 3.0 | $\mathrm{kg\ m^{-1}\ s^{-2}}$ |
| | $\mu$ | 3.0 | $\mathrm{kg\ m^{-1}\ s^{-2}}$ |
| **Other parameters** | | | |
| | $(x_R, y_R, z_R)$ | $(0.75, 0.5, 0.65)$ | $(\mathrm{m}, \mathrm{m}, \mathrm{m})$ |
| | $\Omega$ | $[0, 1.0]^3$ | $\mathrm{m}^3$ |
| | $t_{\mathrm{final}}$ | 0.25 | s |
| | $\Delta t$ | $5.0 \cdot 10^{-4}$ | s |
| | $N_t$ | 500 | – |
| | $N_x$ | 877 945 | – |

# 5 Conclusions

We have presented the MPSA-Newmark method with absorbing boundary conditions for the elastic wave equation. The presentation includes verification and analysis of the MPSA-Newmark method, as well as a demonstration of the method's capabilities through simulation examples.

Convergence analysis of the MPSA-Newmark method was first performed and presented for the case with Dirichlet conditions, where the convergence analysis provided the expected rates of 2 for displacements and between 1 and 2 for tractions. For the case with absorbing boundaries, we presented both convergence and stability analyses of the method. The convergence analyses show that the solution converges with up to order 2 for both displacements and tractions in homogeneous, isotropic, heterogeneous and anisotropic media. The analysis of the stability of the scheme was presented for setups with absorbing boundary conditions on all domain boundaries. We initialized the system with waves that were to propagate out of the domain, and the results showed that energy consistently decreased during time-stepping for all the different waves tested. This indicates that the boundary conditions absorb a significant amount of the waves that hit the boundaries.

Finally, we presented a demonstration of the method's capabilities through simulation examples of wave propagation in the vicinity of a fracture, heterogeneity and anisotropy. The simulation examples show wave reflection in heterogeneous and fractured media, as well as different wave propagation patterns in different directions in anisotropic media. We have thus shown that the present finite volume method exhibits the expected behavior when solving elastic wave propagation problems in anisotropic, heterogeneous and fractured porous media. This indicates that there is potential for a unified finite volume method in poromechanical modelling with the inclusion of seismic effects.

## Acknowledgements:

This project has received funding from the European Research Council (ERC) under the European Union's Horizon 2020 research and innovation program (grant agreement No. 101002507).

## References

Berre, I., Boon, W. M., Flemisch, B., Fumagalli, A., Gläser, D., Keilegavlen, E., Scotti, A., Stefansson, I., Tatomir, A., Brenner, K., Burbulla, S., Devloo, P., Duran, O., Favino, M., Hennicker, J., Lee, I. H., Lipnikov, K., Masson, R., Mosthaf, K.,…Zulian, P. (2021). Verification benchmarks for single-phase flow in three-dimensional fractured porous media. *Advances in Water Resources*, *147*, 103759. https://doi.org/10.1016/j.advwatres.2020.103759

Cardiff, P., & Demirdžić, I. (2021). Thirty Years of the Finite Volume Method for Solid Mechanics. *Archives of Computational Methods in Engineering*, *28*(5), 3721-3780. https://doi.org/10.1007/s11831-020-09523-0

Chopra, A. K. (2012). *Dynamics of Structures*. Pearson Education. https://books.google.no/books?id=eRcvAAAAQBAJ

Clayton, R., & Engquist, B. (1977). Absorbing boundary conditions for acoustic and elastic wave equations. *Bulletin of the Seismological Society of America*, *67*(6), 1529-1540. https://doi.org/10.1785/bssa0670061529

Dalcin, L. D., Paz, R. R., Kler, P. A., & Cosimo, A. (2011). Parallel distributed computing using Python. *Advances in Water Resources*, *34*(9), 1124-1139. https://doi.org/10.1016/j.advwatres.2011.04.013

Demirdzic, I., Martinovic, D., & Ivankovic, A. (1988). Numerical simulation of thermal deformation in welded workpiece. *31*, 209-219.

Dormy, E., & Tarantola, A. (1995). Numerical simulation of elastic wave propagation using a finite volume method. *Journal of Geophysical Research: Solid Earth*, *100*(B2), 2123-2133. https://doi.org/10.1029/94JB02648

Dumbser, M., Käser, M., & De La Puente, J. (2007). Arbitrary high-order finite volume schemes for seismic wave propagation on unstructured meshes in 2D and 3D. *Geophysical Journal International*, *171*(2), 665-694. https://doi.org/10.1111/j.1365-246X.2007.03421.x

Ferguson, J. C., Panerai, F., Borner, A., & Mansour, N. N. (2018). PuMA: the Porous Microstructure Analysis software. *SoftwareX*, *7*, 81-87. https://doi.org/10.1016/j.softx.2018.03.001

Higdon, R. L. (1991). Absorbing boundary conditions for elastic waves. *GEOPHYSICS*, *56*(2), 231-241. https://doi.org/10.1190/1.1443035

Igel, H. (2016). *Computational Seismology: A Practical Introduction*. Oxford University Press.

Jacobsen, I. K., Berre, I., Nordbotten, J. M., & Stefansson, I. (2025). Source Code: A Finite Volume Method for Elastic Waves in Heterogeneous Anisotropic and Fractured Media. https://doi.org/10.5281/zenodo.15056461

Keilegavlen, E., Berge, R., Fumagalli, A., Starnoni, M., Stefansson, I., Varela, J., & Berre, I. (2021). PorePy: an open-source software for simulation of multiphysics processes in fractured porous media. *Computational Geosciences*, *25*(1), 243-265. https://doi.org/10.1007/s10596-020-10002-5

Keilegavlen, E., & Nordbotten, J. M. (2017). Finite volume methods for elasticity with weak symmetry. *International Journal for Numerical Methods in Engineering*, *112*(8), 939-962. https://doi.org/10.1002/nme.5538

Kim, J., Tchelepi, H. A., & Juanes, R. (2011). Stability and convergence of sequential methods for coupled flow and geomechanics: Fixed-stress and fixed-strain splits. *Computer Methods in Applied Mechanics and Engineering*, *200*(13), 1591-1606. https://doi.org/10.1016/j.cma.2010.12.022

Knabner, P., & Angermann, L. (2003). *Numerical Methods for Elliptic and Parabolic Partial Differential Equations*. Springer. https://doi.org/10.1007/b97419

Lemoine, G. I., & Ou, M. Y. (2014). Finite Volume Modeling of Poroelastic-Fluid Wave Propagation with Mapped Grids. *SIAM Journal on Scientific Computing*, *36*(3), B396-B426. https://doi.org/10.1137/130920824

Lemoine, G. I., Ou, M. Y., & LeVeque, R. J. (2013). High-Resolution Finite Volume Modeling of Wave Propagation in Orthotropic Poroelastic Media. *SIAM Journal on Scientific Computing*, *35*(1), B176-B206. https://doi.org/10.1137/120878720

Lie, K. A. (2019). *An Introduction to Reservoir Simulation Using MATLAB/GNU Octave: User Guide for the MATLAB Reservoir Simulation Toolbox (MRST)*. Cambridge University Press. https://doi.org/10.1017/9781108591416

Newmark, N. M. (1959). A Method of Computation for Structural Dynamics. *Journal of the Engineering Mechanics Division*, *85*(3), 67-94. https://doi.org/10.1061/JMCEA3.0000098

Nordbotten, J. M. (2014). Cell-centered finite volume discretizations for deformable porous media. *International Journal for Numerical Methods in Engineering*, *100*(6), 399-418. https://doi.org/10.1002/nme.4734

Nordbotten, J. M. (2016). Stable Cell-Centered Finite Volume Discretization for Biot Equations. *SIAM Journal on Numerical Analysis*, *54*(2), 942-968. https://doi.org/10.1137/15m1014280

Nordbotten, J. M., & Keilegavlen, E. (2021). An Introduction to Multi-Point Flux (MPFA) and Stress (MPSA) Finite Volume Methods for Thermo-Poroelasticity. In D. A. Di Pietro, L. Formaggia, & R. Masson (Eds.), *Polyhedral Methods in Geosciences* (pp. 119-158). Springer International Publishing. https://doi.org/10.1007/978-3-030-69363-3_4

Novikov, A., Voskov, D., Khait, M., Hajibeygi, H., & Jansen, J. D. (2022). A scalable collocated finite volume scheme for simulation of induced fault slip. *Journal of Computational Physics*, *469*, 111598. https://doi.org/10.1016/j.jcp.2022.111598

Payton, R. (2012). *Elastic wave propagation in transversely isotropic media* (Vol. 4). Springer Science & Business Media. https://doi.org/10.1007/978-94-009-6866-0

Seriani, G., & Oliveira, S. P. (2020). Numerical modeling of mechanical wave propagation. *La Rivista del Nuovo Cimento*, *43*(9), 459-514. https://doi.org/10.1007/s40766-020-00009-0

Subbaraj, K., & Dokainish, M. A. (1989). A survey of direct time-integration methods in computational structural dynamics—II. Implicit methods. *Computers & Structures*, *32*(6), 1387-1401. https://doi.org/10.1016/0045-7949(89)90315-5

Tadi, M. (2004). Finite Volume Method for 2D Elastic Wave Propagation. *Bulletin of the Seismological Society of America*, *94*(4), 1500-1509. https://doi.org/10.1785/012003138

Terekhov, K. M. (2020). Cell-centered finite-volume method for heterogeneous anisotropic poromechanics problem. *Journal of Computational and Applied Mathematics*, *365*, 112357. https://doi.org/10.1016/j.cam.2019.112357

Terekhov, K. M., & Tchelepi, H. A. (2020). Cell-centered finite-volume method for elastic deformation of heterogeneous media with full-tensor properties. *Journal of Computational and Applied Mathematics*, *364*, 112331. https://doi.org/10.1016/j.cam.2019.06.047

Terekhov, K. M., & Vassilevski, Y. V. (2022). Finite volume method for coupled subsurface flow problems, II: Poroelasticity. *Journal of Computational Physics*, *462*, 111225. https://doi.org/10.1016/j.jcp.2022.111225

Tsogka, C. (1999). *Modélisation mathématique et numérique de la propagation des ondes élastiques tridimensionnelles dans des milieux fissurés*

Tuković, Ž., Ivanković, A., & Karač, A. (2013). Finite-volume stress analysis in multi-material linear elastic body. *International Journal for Numerical Methods in Engineering*, *93*(4), 400-419. https://doi.org/10.1002/nme.4390

Verwer, J. G., & Sanz-Serna, J. M. (1984). *Convergence of method of lines approximations to partial differential equations*. CWI. Department of Numerical Mathematics [NM].

Virieux, J., Calandra, H., & Plessix, R.-É. (2011). A review of the spectral, pseudo-spectral, finite-difference and finite-element modelling techniques for geophysical imaging. *Geophysical Prospecting*, *59*(5), 794-813. https://doi.org/10.1111/j.1365-2478.2011.00967.x

Zhang, W., Zhuang, Y., & Chung, E. T. (2016). A new spectral finite volume method for elastic wave modelling on unstructured meshes. *Geophysical Journal International*, *206*(1), 292-307. https://doi.org/10.1093/gji/ggw148

Zhang, W., Zhuang, Y., & Zhang, L. (2017). A new high-order finite volume method for 3D elastic wave simulation on unstructured meshes. *Journal of Computational Physics*, *340*, 534-555. https://doi.org/10.1016/j.jcp.2017.03.050

Aavatsmark, I. (2002). An Introduction to Multipoint Flux Approximations for Quadrilateral Grids. *Computational Geosciences*, *6*(3), 405-432. https://doi.org/10.1023/A:1021291114475